%% file: coxrings_v6.tex
\documentclass[12pt]{article}
\usepackage{mathtools}
\usepackage{algorithm}
\usepackage{algorithmic}
\usepackage{mathtools}
\DeclarePairedDelimiter{\ceil}{\lceil}{\rceil}
\usepackage{stmaryrd}
\usepackage[margin=1.05in]{geometry}
\usepackage{titling}
\usepackage{amsthm}
\theoremstyle{plain}
\newtheorem{theorem}{Theorem}[section]
\newtheorem{lemma}{Lemma}[section]
\newtheorem{proposition}{Proposition}[section]
\newtheorem{corollary}{Corollary}[section]
\theoremstyle{definition}
\numberwithin{equation}{section}
\newtheorem{definition}{Definition}[section]

\usepackage{multicol}

\newtheorem{example}{Example}[section]
\newtheorem{remark}{Remark}[section]
\newtheorem{assumption}{Assumption}
\newtheorem{conjecture}{Conjecture}
\newcommand{\N}{\mathbb{N}}
\newcommand{\Conv}{\textup{Conv}}
\newcommand{\Pic}{\textup{Pic}}

\newcommand{\HF}{\textup{HF}}
\newcommand{\Reg}{\textup{Reg}}
\newcommand{\id}{\textup{id}}

\newcommand{\im}{\textup{im}}
\newcommand{\m}{{M}}

\newcommand{\B}{\mathcal{B}}

\newcommand{\Res}{\textup{Res}}
\newcommand{\diag}{\textup{diag}}
\newcommand{\Hom}{\textup{Hom}}
\newcommand{\F}{\mathscr{F}}
\newcommand{\G}{\mathscr{G}}

\newcommand{\Div}{\textup{Div}}
\renewcommand{\div}{\textup{div}}
\newcommand{\codim}{\textup{codim}}
\newcommand{\Spec}{\textup{MaxSpec}}
\newcommand{\Cl}{\textup{Cl}}
\newcommand{\tol}{\textup{tol}}

\newcommand{\A}{\mathscr{A}}

\newcommand{\M}{\mathcal{M}}

\newcommand{\I}{\mathcal{J}}
\newcommand{\K}{K}

\usepackage{url}
\newcommand{\OO}{\mathscr{O}}
\renewcommand{\H}{\mathscr{H}}

\newcommand{\Z}{\mathbb{Z}}

\newcommand{\C}{\mathbb{C}}
\newcommand{\f}{\hat{f}}

\newcommand{\PP}{\mathbb{P}}

\newcommand{\V}{V}
\newcommand{\mean}{\textup{mean}}
\renewcommand{\l}{\ell}
\newcommand{\II}{I}

\newcommand{\z}{\zeta}
\renewcommand{\r}{\rho}
\newcommand{\D}{\delta}
\renewcommand{\z}{\zeta}

\newcommand{\R}{\mathbb{R}}
\usepackage{tikz}
\usepackage{pgfplots}
\usepackage{tikz-3dplot}
\usepackage{mathrsfs}
\newcommand{\Vol}{\textup{Vol}}

\newcommand{\MV}{\textup{MV}}
\providecommand{\keywords}[1]{\small \textbf{Key words ---} #1}
\providecommand{\classification}[1]{\small \textbf{AMS subject classifications ---} #1}
\usepackage{blkarray}
\usepackage{tikz-cd}
\usepackage{multirow}
\usetikzlibrary{external}
\pgfplotsset{
  log x ticks with fixed point/.style={
      xticklabel={
        \pgfkeys{/pgf/fpu=true}
        \pgfmathparse{exp(\tick)}%
        \pgfmathprintnumber[fixed relative, precision=3]{\pgfmathresult}
        \pgfkeys{/pgf/fpu=false}
      }
  }}
\usepackage{array}
\usepackage{enumitem}
\usepackage{authblk}
\usepackage{blindtext}

\usepackage{adjustbox}
\newcolumntype{R}[2]{%
    >{\adjustbox{angle=#1,lap=\width-(#2)}\bgroup}%
    l%
    <{\egroup}%
}
\newcommand{\rot}{\multicolumn{1}{R{45}{1em}}}% no optional argument here, please!
\makeatother
\usepackage{makecell}
\usepackage{amssymb}
\usepackage{subfig}
\makeatletter
\newcommand{\xdashrightarrow}[2][]{\ext@arrow 0359\rightarrowfill@@{#1}{#2}}
\makeatother

\newcommand{\pair}[1]{\langle{#1}\rangle}
\newcommand{\ideal}[1]{\left \langle {#1}  \right \rangle}

\begin{document}
\title{Numerical Root Finding via Cox Rings}
\author[1]{Simon Telen \thanks{\texttt{simon.telen@kuleuven.be}}}
\affil[1]{Department of Computer Science, KU Leuven}
\maketitle

\begin{abstract}
We present a new eigenvalue method for solving a system of Laurent polynomial equations defining a zero-dimensional reduced subscheme of a toric compactification $X$ of $(\C \setminus \{0\})^n$. We homogenize the input equations to obtain a homogeneous ideal $I$ in the Cox ring of $X$ and generalize the eigenvalue, eigenvector theorem for root finding in affine space to compute homogeneous coordinates of the solutions. %We present an algorithm for computing the appropriate multiplication matrices and show how to compute homogeneous coordinates of the solutions from these matrices. 
Several numerical experiments show the effectiveness of the resulting method. In particular, the method outperforms existing solvers in the case of (nearly) degenerate systems with solutions on or near the torus invariant prime divisors. 
\end{abstract}


% REQUIRED
%\begin{AMS}
%14M25, %toric varieties
% 65H04, %roots of polynomial eqns
% 65H10, %systems of equations
% 65H17 %eigenvalues, -vectors
%\end{AMS}

\classification{\small 14M25, %toric varieties
65H04, %roots of polynomial eqns
65H10, %systems of equations
65H17 %eigenvalues, -vectors
}

\section{Introduction}
Many problems in science and engineering can be solved by finding the solutions of a system of (Laurent) polynomial equations. Here, we consider the important case where the number of solutions to the system is finite. There exist many different approaches to tackle this problem \cite{sturmfels2,elkadi_introduction_2007,cattani2005solving}. Symbolic tools such as Groebner bases focus on systems with coefficients in $\mathbb{Q}$ or in finite fields \cite{cox1,sturmfels1996grobner}. For many applications, it is natural to work in finite precision, floating point arithmetic. This is the case, for instance, when the coefficients are known approximately (e.g.\ from measurements) or when it is sufficient to compute solutions accurately up to a certain number of significant decimal digits.
%When, for instance, the involved polynomials are approximations of transcendental nonlinear functions or their coefficients come from measurements, it is natural to work over $\R$ or $\C$. Also, in this type of situations the coefficients are often not known exactly (see Stetter's notion of \textit{empirical polynomials} \cite{stetter}). This motivates the use of finite precision, floating point arithmetic for solving systems of polynomial equations. Another motivation is the fact that using finite precision arithmetic can cause a huge speed-up with respect to symbolic manipulations. Moreover, in many applications it is sufficient to know the solutions up to a certain number of significant decimal digits. 
The most important classes of numerical solvers are homotopy algorithms \cite{bates2013numerically,verschelde1999algorithm,li1997numerical} and algebraic methods such as resultant based algorithms \cite{cox2,emir1,emiris_matrices_1999,noferini,mvb} and normal form algorithms \cite{mourrain1999new,mourrain_stable_2008,dreesen2012back,telen2017stabilized,telen2017solving} which rewrite the problem as an eigenvalue problem. Homotopy solvers are very successful for systems with many variables of low degree, whereas algebraic solvers can handle high degree systems in few variables. The algorithm presented in this paper is a new, numerical normal form algorithm for solving square systems of Laurent polynomial equations. The approach distinguishes itself from existing methods by the interpretation of `solving' the system: we compute the points defined by the input equations on a toric compactification $X$ of $(\C \setminus \{0\})^n \simeq T_X \subset X$ via an eigenvalue computation. More specifically, we work in the Cox ring of $X$ to find `homogeneous' coordinates of the solutions. The motivation is that, even though generically all solutions lie in $(\C \setminus \{0\})^n$, many problems encountered in applications are non-generic with respect to the Newton polytopes of the input equations. Solutions on or near $X \setminus T_X$ cause trouble for the stability of existing numerical algorithms, as we will show in our experiments, and the proposed algorithm is designed to handle such situations. The correctness of the algorithm depends on a conjecture regarding the regularity of a homogeneous ideal in the Cox ring of $X$. In the remainder of this section, we discuss some applications and give an overview of related work and of our main contributions. We conclude the section with an outline of this paper.
%Roughly speaking, path tracking in homotopy continuation employs techniques from numerical ODE solving, whereas algebraic methods rely on numerical linear algebra. 

%A drawback of homotopy continuation is the possibility of path failure due to ill conditioning, path jumping or seemingly diverging paths. See for instance the numerical experiments in \cite{telen2017solving}. These methods give numerical approximations of \textit{all} solutions, given that the system satisfies the necessary genericity assumptions, but they suffer from the curse of dimensionality. 

%The method presented in this paper is a new eigenvalue-based or algebraic method. Therefore, the main applications we have in mind are problems that can be formulated as polynomial equations in only a few variables. 
\subsection*{Applications}
The applications we have in mind are problems that can be formulated as polynomial systems in only a few variables.  

Many problems in computer vision, such as relative pose problems, require the solution of a system of polynomial equations \cite{kukelova2013algebraic,kukelova2008automatic}. In this context, there are often several different polynomial formulations for the same problem, with a different number of variables and a different degree of the equations. See \cite[Sec.\ 7.1.3]{kukelova2013algebraic} for a description of a relative pose problem by a square 7-dimensional system (6 quadratics and a cubic in 7 unknowns) and by a square 3-dimensional system (two cubics and a quintic in 3 unknowns). 

Another application comes from molecular biology. In \cite{emiris1998computer} the problem of computing all possible conformations of several molecules is written in the form of a polynomial system in only two or three variables. 

A problem encountered in many fields of engineering is that of finding the critical points of a function $f$, not necessarily polynomial, in a bounded domain $\Omega \subset \R^n$. A possible approach is to replace $f$ by a polynomial $\tilde{f}$, computed from samples, which approximates $f$ on $\Omega$ and compute the critical points of $\tilde{f}$ instead. The problem is now reduced to a system of polynomial equations, and if $\tilde{f}$ is a good approximation of $f$ in $\Omega$, the solutions in $\Omega$ will be good approximations of the critical points of $f$. It is clear that high degrees lead to better approximations, but also to higher degree polynomial systems. See \cite{noferini} for an application of this technique to solve one of the SIAM 100-Digit Challenge problems \cite{trefethen2002100}.

\subsection*{Related work}
As stated above, solutions on or near the torus invariant prime divisors (i.e.\ the irreducible components of $X \setminus T_X$) cause trouble for numerical root finding in non-compact solution spaces such as $\C^n$ or $(\C \setminus \{0\})^n$. In practice, for homotopy methods, such solutions are the reason for diverging paths, which often require a lot of unnecessary computational effort. Algebraic solvers such as the algorithms proposed in \cite{telen2017stabilized} and \cite[\S 3, \S 4]{telen2017solving}, as well as the classical resultant algorithms \cite[Chapters 3 and 7]{cox2} for computing multiplication matrices, require invertibility of a certain matrix: see for instance the matrix $M_{11}$ in \cite[Chapter 3, \S 6]{cox2} or the matrix $N_{|B}$ in \cite[Section 2]{telen2017solving}. In the presence of solutions on special divisors `at infinity', these matrices are singular. In a numerical context, if these solutions are not exactly \textit{on}, but \textit{near} $X \setminus T_X$, homotopy paths `diverge' to large solutions, causing scaling and condition problems, and the algebraic algorithms require the inversion of an ill-conditioned matrix, causing large rounding errors. A partial solution is to homogenize the equations and solve the problem in $X = \PP^{n_1} \times \cdots \times \PP^{n_k}, k \geq 1, n_1 + \ldots + n_k = n$, which should be thought of as a compactification of $\C^n$, such that a `solution' is defined by $n + k$ (multi-)homogeneous coordinates. This technique is used in total degree homotopies \cite{bates2013numerically,wampler2011numerical}, multihomogeneous homotopies \cite[Chapter 8]{sommese} and in normal form methods such as \cite[\S 5, \S 6]{telen2017solving} or \cite{bender2018towards}. However, depending on the support of the input equations, this standard way of homogenizing may introduce highly singular solutions on the torus invariant divisors, or even destroy 0-dimensionality. More general sparsity structures are taken into account by polyhedral homotopies \cite{li1999solving,hustu,verschelde1994homotopies}, toric or sparse resultants \cite{emir1,cox2,d2002macaulay,pedersen1996mixed,emiris1994monomial,massri2016solving} and truncated normal forms \cite[\S 4]{telen2017solving}. In \cite{huber1998polyhedral} a method for dealing with diverging paths in a polyhedral homotopy is proposed.

In symbolic computing, modified sparse resultant methods have been introduced for solving degenerate systems symbolically \cite{rojas1999solving,d2001computing}. Recently, specialized Groebner basis methods over semigroup algebras have been developed for exploiting sparsity structure \cite{bender2019gr}.  

\subsection*{Contributions}
To the best of the author's knowledge, Cox rings (other than the familiar ones corresponding to products of projective spaces) have not been applied for numerical root finding before. To do so may seem like a bad idea, because the dimension of the Cox ring is (possibly much) greater than that of $X$. However, because of its fine grading by the class group $\Cl(X) = \Div(X)/\sim$ of Weil divisors modulo linear equivalence, this does not affect the computational complexity that much (see Remark \ref{rem:complexity}). The input Laurent polynomial equations define a homogeneous ideal $I$ of the Cox ring $S = \bigoplus_{ \alpha \in \Cl(X) } S_\alpha$ with respect to this grading (this is detailed in Section \ref{sec:setup}). We will assume that $I$ defines a zero-dimensional reduced subscheme $\V_X(I) $ of $X$ which is contained in its largest simplicial open subset $U$ (see Section \ref{sec:preliminaries}). The \textit{regularity} $\Reg(I) \subset \Cl(X)$ of this ideal is defined in Section \ref{sec:lagreg}. In the same section, we conjecture a degree $\alpha \in \Cl(X)$ that is in $\Reg(I)$ (Conjecture \ref{conj}). The correctness of the algorithm depends upon this conjecture, which is supported by some weaker results in Section \ref{sec:lagreg} and by experimental evidence in Section \ref{sec:examples}. For this degree $\alpha \in \Reg(I)$, let $(S/I)_\alpha$ be the degree $\alpha$ part of the graded $S$-module $S/I$. We will construct a linear \textit{multiplication map} $\m_{f} : (S/I)_\alpha \rightarrow (S/I)_\alpha$ with respect to a rational function $f$ on $X$ which is regular at the roots of $I$. Here is a simplified version of Theorem \ref{thm:multiplication}.
\begin{theorem}
Let $V_X(I) = \{\z_1, \ldots, \z_\D\} \subset U$ be reduced and let $\alpha, \alpha_0 \in \Cl(X)$ be such that $\alpha, \alpha+ \alpha_0 \in \Reg(I)$ and there exists $h_0 \in S_{\alpha_0}$ such that $\z_j \notin \V_X(h_0), j = 1, \ldots, \D$. Then for any $g \in S_{\alpha_0}$, the multiplication map $\m_f : (S/I)_\alpha \rightarrow (S/I)_{\alpha}$ with $f = g/h_0$ has eigenvalues $f(\z_j)$.
\end{theorem}
For every monomial $x^{b_i} \in S_{\alpha_0}$, we compute a multiplication matrix and denote its eigenvalues by $\lambda_{ij}, j = 1, \ldots, \D$. This way, we reduce the problem of finding Cox coordinates of $\z_j$ to finding one point on the affine variety defined by the simple binomial system $\{x^{b_i} = \lambda_{ij} ~|~ x^{b_i} \in S_{\alpha_0} \}$ (Corollary \ref{cor:orbiteq}). This leads to a numerical linear algebra based algorithm for finding Cox coordinates (Algorithm \ref{alg:coxcoords}). Unlike other numerical methods, the algorithm is robust in the situation where some of the $\z_j$ are on or near torus invariant prime divisors. We illustrate this in Section \ref{sec:examples} with some examples.

\subsection*{Outline of the paper}
The paper is organized as follows. In the next section we discuss some preliminaries on Cox rings and the classical eigenvalue, eigenvector theorem for polynomial root finding. Our problem setup is discussed in detail in Section \ref{sec:setup}. In Section \ref{sec:lagreg} we introduce homogeneous Lagrange polynomials and their relation to multigraded regularity. Our main result is discussed in detail in Section \ref{sec:toriceval}. The resulting algorithm is presented in Section \ref{sec:alg}. Finally, in Section \ref{sec:examples} we work out several numerical examples. Throughout the paper, we work with polynomials, varieties and vector spaces over $\C$. 

\section{Preliminaries}
\label{sec:preliminaries}
In this section we give a brief introduction to the classical eigenvalue, eigenvector theorem and to complete toric varieties and their Cox rings. We denote by $\V(I) \subset \C^n$ the affine variety of an ideal $I \subset \C[x_1, \ldots, x_n]$ and by $I(Y) \subset \C[x_1, \ldots, x_n]$ the vanishing ideal of a set $Y \subset \C^n$. If $I$ is generated by $f_1, \ldots, f_s \in \C[x_1, \ldots, x_n]$, we denote $I = \ideal{f_1, \ldots, f_s}$ and $V(I) = \V(\ideal{f_1, \ldots, f_s}) = \V(f_1, \ldots, f_s)$. For a finite dimensional vector space $W$, $W^\vee = \Hom_\C(W, \C)$ denotes its dual. For a linear endomorphism $M: W \rightarrow W$ of a finite dimensional vector space $W$, a right eigenpair is $(\lambda, w) \in \C \times (W \setminus \{0\})$ satisfying $M (w) = \lambda w$. Analogously, a left eigenpair is given by $(v, \lambda) \in (W^\vee \setminus \{0 \}) \times \C$  satisfying $v \circ M = \lambda v$ .

\subsection{The classical eigenvalue, eigenvector theorem for polynomial root finding} \label{subsec:multclass}
Let $R = \C[x_1, \ldots, x_n]$ be the ring of $n$-variate polynomials with coefficients in $\C$. Take $f_i \in R, i = 1, \ldots, s$ and let $I = \ideal{f_1, \ldots, f_s}$ be a zero-dimensional ideal in $R$. That is, $\V(I) = \{z_1, \ldots, z_\D\}$ consists of $\D < \infty$ points in $\C^n$. We assume for simplicity that all of the $z_i$ have multiplicity one or, equivalently, that $I$ is radical. By \cite[Chapter 2, Lemma 2.9]{cox2} there exist polynomials $\l_i \in R, i = 1, \ldots, \D$ such that 
$$ \l_i(z_j) = \begin{cases}
0 & i \neq j \\
1 & i = j 
\end{cases}.
$$
The $\l_i$ are called \textit{Lagrange polynomials} with respect to the set $\V(I)$. We define $v_j \in (R/I)^\vee$ by $v_j (f + I) = f(z_j)$. 
\begin{lemma} \label{lem:aff}
The map 
$ \psi : R/I \rightarrow \C^\D : f + I \mapsto (v_1(f+I), \ldots, v_\D(f+I))$
is an isomorphism of vector spaces.
\end{lemma}
\begin{proof}
The map $\psi$ is clearly linear and injective. Surjectivity follows from $\psi(\l_j + I) = e_j$ with $e_j$ the $j$-th standard basis vector of $\C^\D$. 
\end{proof}
It follows from Lemma \ref{lem:aff} that, under our assumptions, $\dim_\C(R/I) = \D$. This is well known, see for instance \cite[Chapter 5, \S3, Proposition 7]{cox1}. In particular, the map $\psi$ defines coordinates on $R/I$ and the residue classes of the Lagrange polynomials form a basis of $R/I$ with dual basis $v_j, j = 1, \ldots, \D$. For $g \in R$, define the linear map $\m_g: R/I \rightarrow R/I: f+I \mapsto fg +I$.
\begin{theorem}[Eigenvalue, eigenvector theorem] \label{thm:EVaff}
The left and right eigenpairs of $\m_g$ are 
$$ (v_j, g(z_j)), \qquad (g(z_j), \l_j + I), \qquad j= 1, \ldots, \D.$$
\end{theorem}
\begin{proof}
See for instance \cite[Chapter 2, Proposition 4.7]{cox2}.
\end{proof}
Note that by definition, $\m_{g_1} \circ \m_{g_2} = \m_{g_2} \circ \m_{g_1}$ for any $g_1, g_2 \in R$. Therefore, after fixing a basis for $R/I$, the matrices corresponding to any two multiplication maps commute and have common eigenspaces. Theorem \ref{thm:EVaff} provides the following algorithm for finding the points in $\V(I)$: 
\begin{enumerate}
\item compute the matrices $\m_{x_1}, \ldots, \m_{x_n}$, 
\item find the coordinates of the $z_i$ from their simultaneous eigenvalue decomposition.
\end{enumerate}
For a more detailed exposition on multiplication matrices, we refer the reader to \cite[Chapter 2]{cox2}, \cite[Chapter 4]{elkadi_introduction_2007} and \cite[Chapter 2]{sturmfels2}.
\subsection{Complete toric varieties and Cox rings} \label{subsec:torvar}
We will restrict ourselves to the discussion of only those aspects of toric varieties that are directly related to this paper. The reader who is unfamiliar with unexplained basic concepts can find an excellent introduction in \cite{cox2011toric} or \cite{fulton1993introduction}. For more information on Cox rings we refer to \cite[Chapter 5]{cox2011toric} and the original paper by Cox \cite{cox1995homogeneous}. The $n$-dimensional algebraic torus $(\C^*)^n = (\C \setminus \{0\})^n$ has character lattice $M = \Hom_\Z((\C^*)^n, \C^*) \simeq \Z^n$ and cocharacter lattice $N = \Hom_\Z(M,\Z) \simeq \Z^n$. An element $m \in M$ gives $\chi^m : (\C^*)^n \rightarrow \C^* $ such that if $m$ corresponds to $(m_1, \ldots, m_n) \in \Z^n$, $\chi^m(t) = t^m = t_1^{m_1} \cdots t_n^{m_n}$. Hence characters can be thought of as Laurent monomials and 
$$\C[M] = \bigoplus_{m \in M} \C \cdot \chi^m \simeq \C[t_1^{\pm 1}, \ldots, t_n^{\pm 1}].$$
Following \cite{cox2011toric}, we denote $N_\R = N \otimes_\Z \R \simeq \R^n$ and $T_N= N \otimes_\Z \C^* = (\C^*)^n $. A complete, normal toric variety $X$ with torus $T_N$ is given by a complete fan $\Sigma$ in $N_\R$ and we will sometimes emphasize this correspondence by writing $X = X_\Sigma$. The set of $d$-dimensional cones of $\Sigma$ is denoted $\Sigma(d)$. In particular, we write $\Sigma(1) = \{\r_1, \ldots, \r_k \}$ for the rays of $\Sigma$ and $u_i \in N$ for the primitive generator of $\rho_i$. It is convenient to think of the $u_i$ as column vectors and to define the matrix $F = [u_1 ~ u_2 ~ \cdots ~ u_k ] \in \Z^{n \times k}$. We will use $F_{ij}$ for the entry in row $i$, column $j$ of $F$, $F_{i,:}$ for the $i$-th row of $F$, $F_{:,j} = u_j$ for the $j$-th column of $F$ and $F^\top$ for the transpose. Every ray $\rho_i$ corresponds to a torus invariant prime divisor $D_i$ on $X_\Sigma$ and we have $X_\Sigma \setminus (\bigcup_{i=1}^k D_i) = T_{X_\Sigma} \simeq T_N$. The class group $\Cl(X_\Sigma)$ of $X_\Sigma$, which is the group of Weil divisors modulo linear equivalence, is generated by the classes $[D_i]$ of the torus invariant prime divisors. The Picard group $\Pic(X_\Sigma) \subset \Cl(X_\Sigma)$ consists of the classes of Weil divisors that are locally principal. Identifying $\bigoplus_{i=1}^k \Z \cdot D_i \simeq \Z^k$ we have a short exact sequence
$$ 0 \longrightarrow M \overset{F^\top}{\longrightarrow} \Z^k \longrightarrow \Cl(X_\Sigma) \longrightarrow 0$$
where $\Z^k \longrightarrow \Cl(X_\Sigma)$ sends a torus invariant Weil divisor $\sum_{i=1}^k a_i D_i$ to its class $ [ \sum_{i=1}^k a_i D_i] \in \Cl(X_\Sigma)$. Taking $\Hom_\Z(-, \C^*)$ and defining the \textit{reductive group} $G = \Hom_\Z(\Cl(X_\Sigma), \C^*)$ we find that $G$ is the kernel of the map 
\begin{equation} \label{eq:GCQ}
 \pi: (\C^*)^k \rightarrow T_N : t \mapsto (t^{F_{1,:}}, \ldots, t^{F_{n,:}}).
 \end{equation}
That is, $G$ is the subgroup of $(\C^*)^k$ given by 
%Each $\sigma \in \Sigma(n)$ defines an affine toric variety $U_\sigma = \Spec(\C[\sigma^\vee \cap M])$. The abstract toric variety $X_\Sigma$ is obtained by gluing together the affine toric varieties $\{ U_{\sigma}, \sigma \in \Sigma \}$ along the open subsets corresponding to common faces \cite{cox2011toric}. The torus $T_{X_\Sigma} \subset X_\Sigma$ is isomorphic to $(\C^*)^n$ and it is Zariski dense in $X_\Sigma$.

%\begin{definition}[Reductive group]
%Define
$$ G = \{ g \in (\C^*)^k : g^{F_{i,:}} = 1, i = 1, \ldots, n \}$$
and $\pi$ is constant on $G$-orbits.
%The group $G$ is called the \textit{reductive group} of $\Sigma$ (or of $X_\Sigma$).
%\end{definition}
Let $S = \C[x_1, \ldots, x_k]$ be the polynomial ring in $k$ variables where each of the $x_i$ corresponds to a ray $\rho_i \in \Sigma(1)$. For every cone $\sigma \in \Sigma$, denote by $\sigma(1)$ the rays contained in $\sigma$. We are going to associate a monomial in $S$ to each cone in $\Sigma$: for $\sigma \in \Sigma$, define 
$ x^{\hat{\sigma}} = \prod_{\rho_i \notin \sigma(1)} x_i.$
The \textit{irrelevant ideal} $\K$ of $\Sigma$ (or of $X_\Sigma$) is the monomial ideal defined as 
\begin{equation} \label{eq:irrelideal}
\K = \ideal{ x^{\hat{\sigma}} : \sigma \in \Sigma(n)} \subset S.
\end{equation}
The \textit{exceptional set} of $X_\Sigma$ is $Z = \V(\K) \subset \C^k$. The action of $G$ on $(\C^*)^k$ extends to an action on $\C^k \setminus Z$. In \cite{cox1995homogeneous}, Cox proves that there is a good categorical quotient 
$ \pi : \C^k \setminus Z \rightarrow X_\Sigma $, constant on $G$-orbits, such that \eqref{eq:GCQ} is its restriction to $(\C^*)^k$.  By the properties of good categorical quotients we have a bijection 
\begin{equation*} \label{eq:bijection}
\{ \textup{ closed $G$-orbits in $\C^k \backslash Z$ } \} \leftrightarrow \{ \textup{ points in $X_\Sigma$ } \}.
\end{equation*}
Moreover, $\pi$ is an almost geometric quotient, meaning that there is a Zariski open subset $U \subset X_\Sigma$ such that $\pi_{|\pi^{-1}(U)}: \pi^{-1}(U) \rightarrow U$ is a geometric quotient:
$$ \{ \textup{ $G$-orbits in $\pi^{-1}(U)$ } \} \leftrightarrow \{ \textup{ points in $U$ } \}.$$
The open set $U$ is the toric variety $X_{\Sigma'} \subset X_\Sigma$ corresponding to the subfan $\Sigma' \subset \Sigma$ of simplicial cones of $\Sigma$ (see \cite[proof of Theorem 5.1.11]{cox2011toric}). Therefore, by the orbit-cone correspondence, $X \setminus U$ is a union of $T_N$-orbits of codimension at least 3 (cones of dimension 0, 1 or 2 are simplicial). If $\Sigma$ is simplicial, the nicest possible bijection holds: 
$$ \{ \textup{ $G$-orbits in $\C^k \backslash Z$ } \} \leftrightarrow \{ \textup{ points in $X_\Sigma$ } \}.$$
In this case we write $X_\Sigma = (\C^k \setminus Z)/G$.
\begin{example}
The quotient construction of $X_\Sigma$ is a generalization of the familiar construction of $\PP^n$ as the quotient $\PP^n = (\C^{n+1} \setminus \{0\})/ \C^*$. In this case $S = \C[x_0, \ldots, x_n]$, $\K = \ideal{x_0, \ldots, x_n}$, $Z = \{0\}$ and $G = \Hom_\Z(\Cl(\PP^n), \C^*) = \Hom_\Z(\Z, \C^*) = \C^*$ acts by $g \cdot (x_0, \ldots, x_n) = (g x_0, \ldots, g x_n), g \in G$.
\end{example} 
The ring $S$ has a natural grading by $\Cl(X_\Sigma)$: 
\begin{equation} \label{eq:grading}
\deg(x^a) = \deg(x_1^{a_1} \cdots x_k^{a_k}) = [ \sum_{i=1}^k a_i D_i ] \in \Cl(X_\Sigma), \qquad S = \bigoplus_{\alpha \in \Cl(X_\Sigma)} S_\alpha,
\end{equation}
where $S_\alpha = \bigoplus_{\deg(x^a) = \alpha} \C \cdot x^a$. In fact, the only nonzero graded pieces correspond to `positive' degrees, and one can write 
$$\Cl(X_\Sigma)_+ = \{ \alpha \in \Cl(X_\Sigma) ~|~ \alpha = n_1 \deg(x_1) + \cdots + n_k \deg(x_k), n_i \in \N \}, \quad S = \bigoplus_{\alpha \in \Cl(X_\Sigma)_+} S_\alpha.$$
Similarly, we denote $\Pic(X_\Sigma)_+ = \Cl(X_\Sigma)_+ \cap \Pic(X_\Sigma)$. The graded pieces correspond to vector spaces of global sections of divisorial sheaves, that is, for $\alpha \in \Cl(X_\Sigma)$ with $\alpha = [ D], D = \sum_{i=1}^k a_i D_i $, 
\begin{equation} \label{eq:gradedpiece}
S_\alpha \simeq \Gamma( X_\Sigma, \OO_{X_\Sigma}(D)) \simeq \bigoplus_{ F^\top m + a \geq 0} \C \cdot \chi^m.
\end{equation}
Here the direct sum ranges over all $m$ such that elementwise, $F^\top m + a \geq 0$, that is, $\pair{u_i,m} + a_i \geq 0, i = 1, \ldots, k$ where $\pair{\cdot, \cdot}$ is the natural pairing between $N$ and $M$. Denoting $x^{F^\top m + a} = x_1^{\pair{u_1,m} + a_1} \cdots x_k^{\pair{u_k,m} + a_k}$, the isomorphism \eqref{eq:gradedpiece} is given by 
\begin{equation} \label{eq:hom}
\sum_{F^\top m + a \geq 0} c_m \chi^m \mapsto \sum_{F^\top m + a \geq 0} c_m x^{F^\top m + a} \in S_{\alpha}, 
\end{equation}
which is \textit{homogenization} with respect to $\alpha$. To see the analogy with the classical notion of homogenization, note that the action of $G$ on $\C^k$ induces an action of $G$ on $S$ by $(g \cdot f )(x) = f(g^{-1} \cdot x)$ for $g \in G, f \in S$. If $f \in S_\alpha$, it is the image of some Laurent polynomial under \eqref{eq:hom} and we can write
\begin{equation} \label{eq:eigenspaces}
(g \cdot f)(x) = \sum_{F^\top m + a \geq 0} c_m (g^{-1} \cdot x)^{F^\top m + a} = g^{-a} f(x)
\end{equation}
since by the definition of the reductive group $g^{F^\top m}= 1$. This shows that the number $g^{-a}$ does not depend on the representative divisor $D$ we choose for $\alpha \in \Cl(X_\Sigma)$. It therefore makes sense to write $g^{-\alpha} = g^{-(F^\top m +a)}$. Equation \eqref{eq:eigenspaces} shows that the homogeneous components $S_\alpha \subset S$ with respect to the grading \eqref{eq:grading} are the eigenspaces of the action of $G$ on $S$ and that
\begin{equation} \label{eq:zeroset}
\V_{X_\Sigma}(f) = \{ p \in X_\Sigma : f(x) = 0 \textup{ for some } x \in \pi^{-1}(p)\} \subset X_\Sigma
\end{equation} 
is well defined if $f$ is homogeneous. An ideal $I \subset S$ is called homogeneous if it is generated by homogeneous polynomials, and it is straightforward to extend \eqref{eq:zeroset} to define $\V_{X_\Sigma}(I)$. The ring $S$ equipped with the grading \eqref{eq:grading} and the irrelevant ideal \eqref{eq:irrelideal} is called the \textit{total coordinate ring}, \textit{homogeneous coordinate ring} or \textit{Cox ring} of $X_\Sigma$. 
\begin{example}
The complete fans $\Sigma$ we will encounter in this paper are normal fans of full dimensional lattice polytopes \cite[\S 2.3]{cox2011toric}. If 
$$ P = \{ m \in M_\R ~|~ \pair{u_i,m} \geq -a_i, i = 1, \ldots, k \}$$
is the minimal facet representation of a full dimensional lattice polytope $P \subset M_\R$, then its normal fan $\Sigma_P$ defines a toric variety $X_{\Sigma_P}$, which we will often denote by $X$ for simplicity of notation. There are bijective correspondences between rays in $\Sigma_P$, facets of $P$, torus invariant prime divisors in $X$ and indeterminates in the Cox ring. The matrix $F$ contains the primitive inward pointing facet normals of $P$. For example, the toric variety of the standard $n$-simplex is $\PP^n$.
\end{example}

\begin{example} \label{ex:hirz1}
As a running example, we will consider the problem of finding the intersections of two curves on the Hirzebruch surface $\H_2$. The associated fan $\Sigma$ and the matrix $F$ of ray generators are shown in Figure \ref{fig:hirzebruch1}. The Cox ring $S = \C[x_1, x_2,x_3,x_4]$ is graded by $ \Cl(\H_2)\simeq \Z^4/\im F^\top \simeq \Z^2$, with $\deg(x^b) = \deg(x_1^{b_1}x_2^{b_2}x_3^{b_3}x_4^{b_4}) = (b_1 - 2b_2+b_3, b_2 + b_4)$. The reductive group and exceptional set are given by 
$G = \{(\lambda, \mu, \lambda, \lambda^2 \mu) ~|~ (\lambda,\mu)\in (\C^*)^2 \} \subset (\C^*)^4$ and $Z = \V(x_1,x_3) \cup \V(x_2,x_4) \subset \C^4$ respectively.
Since $\H_2$ is smooth, it is simplicial (in the notation from above $U = \H_2$) and $\Pic(\H_2) = \Cl(\H_2)$.
\begin{figure}
\centering
\input{fan.tex}
\caption{Fan and matrix of primitive ray generators of the Hirzebruch surface $\H_2$.}
\label{fig:hirzebruch1}
\end{figure}

\end{example}

\section{Problem setup} \label{sec:setup}
In this section, we give a detailed description of the problem considered in this paper and we discuss our assumptions. We start from $n$ given Laurent polynomials $\f_1, \ldots, \f_n \in \C[M]$ (that is, we consider \textit{square} systems). Denote
$$ \f_j = \sum c_{m,j} \chi^m$$
and let $P_j \subset M_\R$ be the Newton polytope of $\f_j$: $P_j = \Conv(m \in M ~|~ c_{m,j} \neq 0) \subset M_\R$. Let $P = P_1 + \ldots + P_n$ be the Minkowski sum of these polytopes. We assume that $P$ is full-dimensional and we let $X = X_{\Sigma_P}$ be the complete normal toric variety corresponding to its normal fan. To each $P_j$, we associate a basepoint free\footnote{For $\alpha = [D] \in \Pic(X)$, we say that $p \in X$ is a \textit{basepoint} of $S_\alpha \simeq \Gamma(X,\OO_X(D))$ if every global section of the associated line bundle $\OO_X(D)$ vanishes at $p$. The divisor $D$ and its associated degree $\alpha \in \Pic(X)$ are called \textit{basepoint free} if $S_\alpha$ has no basepoints.} Cartier divisor $D_{P_j}$ on $X$, given by 
$$ D_{P_j} = \sum_{i=1}^k a_{j,i} D_i, \qquad a_{j,i} = - \min_{m \in P_j} \pair{u_i,m}$$
and we denote $a_j = (a_{j,1}, \ldots, a_{j,k}) \in \Z^k, [D_{P_j}] = \alpha_j \in \Pic(X)$. %By construction, the $D_{P_j}$ are basepoint free divisors on $X$, meaning that the morphisms $X \rightarrow \PP(S_{\alpha_j}^\vee)$ have no basepoints \cite[\S 6.0]{cox2011toric}. The following result will be useful. 
For this construction, $D_{P_i+ P_j} = D_{P_i} +  D_{P_j}$ and for $\I \subset \{1,\ldots, n\}$, $P_\I = \sum_{j \in \I} P_j$ we have 
\begin{equation} \label{eq:sumofdiv}
\Gamma(X, \OO_X( \sum_{j \in \I} D_{P_j})) =\Gamma(X, \OO_X(  D_{P_\I}))  = \bigoplus_{m \in P_\I \cap M } \C \cdot \chi^m.
\end{equation}
%\begin{lemma} \label{lem:sumdivsumpol}
%For the above construction, we have $D_{P_i+ P_j} = D_{P_i} +  D_{P_j}$. 
%\end{lemma}
%\begin{proof}
%Since $-a_{i,\ell} = \min_{m \in P_i} \pair{u_\ell,m}$ (the minimum is attained in $P_i$ by the construction) and analogously for $a_{j,\ell}$, we have by \cite[Theorem 3.1.6]{weibel2007minkowski} that $\min_{m \in P_i + P_j} \pair{u_\ell,m} = \min_{m \in P_i} \pair{u_\ell,m}+\min_{ P_j} \pair{u_\ell,m} = - (a_{i,\ell} + a_{j,\ell})$. Since the normal fan $\Sigma_{P_i + P_j}$ of $P_i + P_j$ is the coarsest common refinement of $\Sigma_{P_i}$ and $\Sigma_{P_j}$ \cite[Proposition 6.2.13]{cox2011toric} and $\Sigma_P$ is a refinement of $\Sigma_{P_i + P_j}$, we have 
%$$ P_i + P_j = \{ m \in M_\R ~|~ F^\top m + (a_i + a_j) \geq 0 \}.$$ 
%It follows directly that $D_{P_i} + D_{P_j} = D_{P_i+P_j}$. 
%\end{proof}
%\begin{corollary}
%Let $\I \subset \{1, \ldots, n \}$ and let $P_\I = \sum_{j \in \I} P_j$. Then
%\begin{equation}
%\dim_\C \left (\Gamma(X, \OO_X( \sum_{j \in \I} D_{P_j}) \right ) = |P_\I \cap M |.
%\end{equation}
%\end{corollary}
%\begin{proof}
%By Lemma \ref{lem:sumdivsumpol}, $\sum_{j \in \I} D_{P_j}= D_{P_\I}$ and by \eqref{eq:gradedpiece}
%$$ \Gamma(X, \OO_X(D_{P_\I})) \simeq \bigoplus_{F^\top m + a_\I \geq 0} \C \cdot \chi^m = \bigoplus_{m \in P_\I \cap M} \C \cdot \chi^m,$$
%with $a_\I = \sum_{j \in \I} a_j$. 
%\end{proof}
By definition, $m \in P_j \cap M$ if and only if $F^\top m + a_j \geq 0$, so we have
\begin{equation} \label{eq:linebundle}
\f_j = \sum_{m \in P_j \cap M} c_{m,j} \chi^m \in \Gamma(X, \OO_X(D_{P_j})).
\end{equation}
Homogenizing with respect to $\alpha_j$ according to \eqref{eq:hom} gives (see \cite{cattani1997global})
$$\f_j \mapsto f_j = \sum_{m \in P_j \cap M} c_{m,j} x^{F^\top m + a_j} \in S_{\alpha_j}.$$
Equation \eqref{eq:linebundle} shows that $\f_j$ is a global section of the line bundle given by $\OO_X(D_{P_j})$ \cite[Chapter 6]{cox2011toric}. Its divisor of zeroes is the effective divisor $\div(\f_j) + D_{P_j}$, whose support is exactly $\V_X(f_j)$. This construction gives a homogeneous ideal $I = \ideal{f_1, \ldots, f_n} \subset S$. We will make the following assumptions on $I$.
\begin{assumption} \label{ass:1}
$\V_X(I)$ is zero-dimensional. We denote $\V_X(I) = \{\z_1, \ldots, \z_\D \} \subset X$. 
\end{assumption}
\begin{assumption} \label{ass:2}
$\V_X(I) \subset U \subset X$, where $U$ is the `simplicial part' of $X$ as in Subsection \ref{subsec:torvar}.
\end{assumption}
\begin{assumption} \label{ass:3}
$I$ defines a reduced subscheme of $U \subset X$. That is, all points $\z_i$ are `simple roots' of $I$.
\end{assumption}
%In summary, we assume the solutions of $I$ on $X$ to be finite in number, we assume that they do not lie on a subvariety of $X$ of codimension at least 3 and we assume them to be nonsingular. 
%The last assumption is made for simplicity of the exposition. The results are expected to generalize to the case of `fat points' with multiplicities just as in the classical case (see for instance \cite[Chapter 4, Proposition 2.7]{cox2}). 
It is clear that when $n = 2$, Assumption \ref{ass:2} can be dropped. For $n = 3$, $U$ is the complement of finitely many points in $X$: one point for each vertex of $P$ corresponding to a non-simplicial, full dimensional cone of $\Sigma_P$. It follows that we can drop Assumption \ref{ass:2} also for $n= 3$, since `face systems' corresponding to vertices do not contribute any solutions (see for instance the appendix in \cite{hustu}).
%Instead of assuming $\V_X(I) \subset U$, we could just replace $X$ by $U$ or assume that $X$ is simplicial%(note that $U$ has the same Cox ring as $X$ and under our assumptions, $\V_U(I) = \V_X(I)$)
For $n>3$, Assumption \ref{ass:2} can be dropped if $X$ is simplicial. %We formulate the assumption as above to emphasize that we do not need the normal fan of $P$ to be simplicial.
We will comment on Assumption \ref{ass:3} in Section \ref{sec:lagreg} (Remark \ref{rem:multiplicities}).

In order to say something more about the number $\D$ in Assumption \ref{ass:1}, we recall the definition of mixed volume.
The $n$-dimensional \textup{mixed volume} of a collection of $n$ polytopes $P_1,\ldots,P_n$ in $M_\R \simeq \R^n$, denoted $\MV(P_1,\ldots,P_n)$, is the coefficient of the monomial $\lambda_1 \lambda_2 \cdots \lambda_n$ in $\Vol_n(\sum_{i = 1}^n \lambda_i P_i)$.
A formula for the mixed volume that will be useful is (see \cite{bihan2016irrational,csahin2016multigraded})
\begin{equation} \label{eq:mvpoints}
\MV(P_1, \ldots,P_n) = \sum_{\ell = 0}^{n} (-1)^{n-\ell} \sum_{ \substack{\I \subset \{1,\ldots,n\} \\ |\I| = \ell}} \left|(P_0+P_\I) \cap M \right |,
\end{equation}
for any lattice polytope $P_0 \subset \R^n$ corresponding to a basepoint free divisor $D_{P_0}$. The following important theorem was named after Bernstein, Khovanskii and Kushnirenko and tells us what the number $\D$ is. 
\begin{theorem}[BKK Theorem] \label{thm:bkk}
Let $I = \ideal{ f_1,\ldots,f_n } \subset S$ be a homogeneous ideal constructed as above. If $I$ defines $\D < \infty$ points on $X$, counting multiplicities, then $\D$ is given by $\MV(P_1,\ldots,P_n)$. For generic choices of the coefficients of the $f_i$, the number of roots in $T_X \simeq T_N  = (\C^*)^n$ is exactly equal to $\MV(P_1,\ldots,P_n)$ and they all have multiplicity one.
\end{theorem}
\begin{proof}
See \cite[\S 5.5]{fulton1993introduction}. For sketches of the proof we refer to \cite{cox2,sturm}. Other proofs can be found in Bernstein's original paper \cite{bernstein} and in \cite{hustu}.
\end{proof}
Theorem \ref{thm:bkk} is a generalization of B\'ezout's theorem for projective space. Motivated by this result, for the rest of this article $\D = \MV(P_1, \ldots, P_n)$.
We can represent each $\z_j \in \V_X(I)$ by a set of homogeneous coordinates $z_j = (z_{j1}, \ldots, z_{jk}) \in \C^k \setminus Z$. Let $\pi^{-1}(\z_j) = G \cdot z_j \subset \C^k \setminus Z$ be the corresponding $(k - n)$-dimensional closed $G$-orbit and let $\overline{G \cdot z_j}$ be the closure in $\C^k$. It follows from our assumptions that 
$$ \V(I)  \setminus Z = G \cdot z_1 \cup \cdots \cup G \cdot z_\D \qquad \textup{and} \qquad \V(I) = \overline{G \cdot z_1} \cup \cdots \cup \overline{G \cdot z_\D} \cup Z',$$
with $Z' \subset Z$ a closed subvariety. We define $J = \II(\overline{G \cdot z_1} \cup \cdots \cup \overline{G \cdot z_\D} )$ to be the ideal of the union of orbit closures, which is radical and saturated with respect to the irrelevant ideal $K$. The ideal $J$ is the one investigated in \cite{csahin2016multigraded} (in the simplicial case). It is clear that $I \subset J$. In some special cases where $Z$ is very small, the ideals $I$ and $J$ coincide. This happens for instance for $X = \PP^n$ or for any weighted projective space $X = \PP(w_0, \ldots, w_n)$.

\begin{example} \label{ex:hirz2}
Let us consider the polynomials
\begin{eqnarray*}
\f_1 &=& 1 + t_1 +t_2 + t_1t_2 + t_1^2t_2 + t_1^3t_2 , \\
\f_2 &=& 1 + t_2 + t_1t_2 + t_1^2t_2.
\end{eqnarray*}
We think of $\f_1,\f_2$ as elements of $\C[t_1^{\pm 1}, t_2^{\pm 1}] \simeq \C[M]$ with $M = \Z^2$ the character lattice of $T_N = (\C^*)^2$. The polytopes $P_1, P_2$ and $P$ are shown in Figure \ref{fig:hirzebruch2}. Note that the normal fan $\Sigma_P$ of $P$ is the fan of Figure \ref{fig:hirzebruch1}, so the toric variety associated to this system is $X = X_{\Sigma_P} = \H_2$.
\begin{figure}
\centering
\input{newtonpolygons.tex}
\caption{Newton polytopes involved in Example \ref{ex:hirz2}.}
\label{fig:hirzebruch2}
\end{figure}
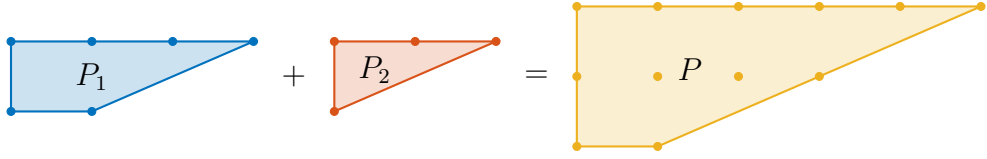
We identify $\Cl(X)$ with $\Z^2$ as in Example \ref{ex:hirz1}. It is easy to check that $\alpha_1 = [D_{P_1}] = [D_3+D_4] = (1,1) \in \Cl(X)$ and $\alpha_2 = [D_{P_2}] = [D_4] = (0,1) \in \Cl(X)$. This gives the following homogeneous polynomials in the Cox ring $S = \C[x_1, \ldots, x_4]$: 
\begin{eqnarray*}
f_1 &=& x_3x_4 + x_1x_4 + x_2x_3^3 + x_1x_2x_3^2 + x_1^2x_2x_3 + x_1^3x_2, \\
f_2 &=& x_4 + x_2x_3^2 + x_1x_2x_3 + x_1^2x_2.
\end{eqnarray*}
The mixed volume is $\D = \MV(P_1,P_2) = 3$. To see that the ideal $I = \ideal{f_1,f_2}$ satisfies our assumptions, we compute its primary decomposition\footnote{We used Macaulay2 to perform the symbolic computations in this example \cite{eisenbud2001computations}.}.
\begin{eqnarray*}
I &=& \ideal{x_1+x_3, x_2x_3^2+x_4}\cap \ideal{x_1,x_2x_3^2+x_4} \cap \ideal{x_3, x_1^2x_2+x_4} \cap \ideal{x_2,x_4}
\end{eqnarray*}
which gives the decomposition of the associated variety $\V(I) = \overline{G \cdot z_1} \cup \overline{G \cdot z_2} \cup \overline{G \cdot z_3} \cup Z'$ with orbit representatives
$ z_1 = (-1,-1,1,1), z_2 = (0,-1,1,1), z_3 = (1,-1,0,1)$ and $Z' = \V(x_2,x_4) \subset Z$. This shows that $I$ defines the expected number of simple, isolated points on $X = \H_2$. The first solution $\z_1 = \pi(z_1) \in T_N$ lies in the torus, the others satisfy $\z_2 = \pi(z_2) \in D_1$, $\z_3 = \pi(z_3) \in D_3$. The ideal $J$ in this example is the intersection of the first three primary components of $I$. We find $J = \ideal{x_1^2x_3+x_1x_3^2, f_2}$.
\end{example}

\section{Multigraded regularity and homogeneous Lagrange polynomials}
\label{sec:lagreg}
The regularity of a graded module measures its complexity (for instance, in terms of the degree of minimal generators). The notion of regularity has been studied in a multigraded context. The general situation is treated in \cite{maclagan2003multigraded}. The zero-dimensional case is further investigated in \cite{csahin2016multigraded} and some more results in a multiprojective setting can be found in \cite{bender2018towards,sidman2006multigraded}. In our case, the regularity (as defined below) of the ideal $I$ in Section \ref{sec:setup} will determine in which graded piece $S_\alpha$ of the Cox ring $S$ we can work to define our multiplication maps in Section \ref{sec:toriceval}. The `larger' this graded piece (i.e.\ the larger the dimension of $S_\alpha$ as a $\C$-vector space), the larger the matrices involved in the presented algorithm in Section \ref{sec:alg}. We will define homogeneous Lagrange polynomials and show how they are related to multigraded regularity. As in Subsection \ref{subsec:multclass}, these Lagrange polynomials and their dual basis will have a nice interpretation as eigenvectors of multiplication maps. For $\alpha \in \Cl(X)$, we denote $n_\alpha = \dim_\C(S_\alpha)$. Since $X$ is complete, $n_{\alpha} < \infty, \forall \alpha \in \Cl(X)$ \cite[Proposition 4.3.8]{cox2011toric}. The ideals $I, J \subset S$ are as defined in Section \ref{sec:setup}. In particular, $I$ satisfies Assumptions \ref{ass:1}-\ref{ass:3}. For $\alpha \in \Cl(X)$, let $S_\alpha = \bigoplus_{i=1}^{n_\alpha} \C \cdot x^{b_i}, b_i \in \N^k$ and consider the map $$\Phi_\alpha : \C^k \setminus Z \dashrightarrow \PP^{n_\alpha - 1} \simeq \PP(S_\alpha^\vee) \simeq \PP(\Gamma(X,\OO_X(D))^\vee):(x_1, \ldots, x_k) \mapsto (x^{b_1}, \ldots, x^{b_{n_\alpha}} ).$$
Note that $\Phi_\alpha$ may have basepoints (hence the dashed arrow) and it is constant on $G$-orbits. We will say that $\alpha \in \Cl(X)$ is basepoint free if $\Phi_\alpha$ has no basepoints (this extends the definition for basepoint free $\alpha \in \Pic(X)$ to the class group). We say that $\z \in U \subset X$ is a basepoint of $S_\alpha$ if $\pi^{-1}(\z)$ are basepoints of $\Phi_\alpha$. %If $\alpha$ is basepoint free, by the universal property of a good categorical quotient the map $\Phi_\alpha$ factors as $\Phi_\alpha = \phi_\alpha \circ \pi$, with $\phi_\alpha : X \rightarrow \PP^{n_{\alpha} - 1}$. 
The following lemma is straightforward and we omit the proof. 
%We assume that $\Phi_\alpha$ is well defined at the $z_j$. We start from the following Lemma. 
\begin{lemma} \label{lem:nonzero}
Let $\alpha = [D] \in \Cl(X)$ be such that no $\z_j$ is a basepoint of $S_\alpha$. For generic $h \in S_\alpha$, we have $\z_j \notin \V_X(h), j = 1, \ldots, \D$.  
\end{lemma}
%\begin{proof}
%The characters $\chi^m \in \Gamma(X,\OO_X(D))$ form a $\C$-basis of $\Gamma(X,\OO_X(D))$. For each $\z_j$, at least one of these sections does not vanish. It follows that a generic linear combination $\hat{h}$ does not vanish at any of the $\z_j$ either. Take $h$ to be the homogenization in degree $\alpha$ of $\hat{h}$.
%\end{proof}
Note that in particular, the condition of Lemma \ref{lem:nonzero} is always satisfied for basepoint free $\alpha$. The grading on $S$ defines a grading on the quotient $S/I$: $(S/I)_\alpha = S_\alpha / I_ \alpha$. It follows from Lemma \ref{lem:nonzero} that for any $\alpha = [D] \in \Cl(X)$ such that no $\z_j$ is a basepoint of $S_\alpha$, the following $\C$-linear map is well defined for generic $h \in S_\alpha$: 
\begin{equation} \label{eq:evalmap}
\psi_\alpha : (S/I)_\alpha \rightarrow \C^\D : f + I_\alpha \mapsto \left ( \frac{f}{h}(z_1), \ldots, \frac{f}{h}(z_\D) \right ).
\end{equation}
%By Lemma \ref{lem:nonzero}, $\psi_\alpha$ is well defined for generic $h \in S_\alpha$ with basepoint free $\alpha \in \Pic(X)$. 
We fix such a generic $h \in S_\alpha$. Note that the definition of $\psi_\alpha$ does not depend on the choice of representative $z_j$ of $G \cdot z_j$. %The function $(f/h)$ is a rational function on $X$
%(it is a morphism $U \setminus \V_X(h) \rightarrow \C$, and $(X \setminus \V_X(h)) \setminus (U \setminus \V_X(h))$ is of codimension $>$ 2 in $X \setminus \V_X(h)$ \cite[Corollary 4.0.15]{cox2011toric})
%We write $\C(X)$ for the field of rational functions on $X$ and we have $(f/h) \in \C(X)$. 
We will now investigate for which $\alpha \in \Cl(X)$ the map $\psi_\alpha$ defines coordinates on $(S/I)_\alpha$, that is, for which $\alpha$ it is an isomorphism (note that this is independent of the choice of $h$ satisfying $\z_j \notin V_X(h)$). It is clear that for this to happen, we need $\dim_\C((S/I)_\alpha) = \D$. The dimension of the graded parts of $S/I$ is given by the multigraded analog of the Hilbert function \cite{csahin2016multigraded}.
\begin{definition}[Hilbert function] \label{def:hilbfun}
For a homogeneous ideal $I$ in the Cox ring $S$ of $X$, the \textup{Hilbert function} of $I$ is given by 
$ \HF_I : \Cl(X) \rightarrow \N : \alpha \mapsto \dim_\C((S/I)_\alpha).$
\end{definition}
We note that in \cite{csahin2016multigraded}, the Hilbert function of the scheme $\V_X(I)$ is equal to $\HF_J$ as defined above. In order to state a necessary and sufficient condition for surjectivity of $\psi_\alpha$, we will introduce a homogeneous analog of the Lagrange polynomials introduced in Subsection \ref{subsec:multclass}.
\begin{definition}[homogeneous Lagrange polynomials] \label{def:lagpol}
Let $\alpha \in \Cl(X)$ be such that no $\z_j$ is a basepoint of $S_\alpha$ and let $h \in S_\alpha$ be such that $\z_j \notin \V_X(h), j = 1, \ldots, \D$. A set of elements $\l_1, \ldots, \l_\D \in S_\alpha$ is called a set of \textup{homogeneous Lagrange polynomials} of degree $\alpha$ with respect to $h$ if for $j = 1, \ldots, \D$, 
\begin{enumerate}
\item $\z_i \in \V_X(\l_j), i \neq j$, 
\item $\z_j \in \V_X(h-\l_j)$.
\end{enumerate}
\end{definition}
In terms of the homogeneous coordinates $z_j$, a set of homogeneous Lagrange polynomials satisfies $\l_j(z_i) = 0, i \neq j$ and $\l_j(z_j) = h(z_j), j = 1, \ldots, \D$. %Note that any set of polynomials satisfying the first condition of Definition \ref{def:lagpol} and $\z_j \notin \V_X(\l_j)$ can be turned into a set of homogeneous Lagrange polynomials by performing an appropriate scaling. 
\begin{remark} \label{rem:multiplicities}
Let $\l_j, j = 1, \ldots, \D$ be a set of homogeneous Lagrange polynomials of degree $\alpha$ with respect to $h$. The cosets $\ell_j + I_\alpha \in (S/I)_\alpha$ are a dual basis for the evaluation functionals $v_j \in (S/I)_\alpha^\vee$ given by $v_j: (S/I)_\alpha \rightarrow \C : f + I_\alpha \mapsto (f/h)(z_j)$. If $I$ defines points with multiplicities (the case of `fat points', violating Assumption \ref{ass:3}), a starting point would be to extend this set of evaluation functionals to a basis of $(S/I)_\alpha^\vee$, using analogs of differentiation operators. It is known that the theory for the affine root finding problem (Subsection \ref{subsec:multclass}) extends nicely in this way; see for instance \cite[Chapter 4, Proposition 2.7]{cox2}, \cite[Section 4.3]{elkadi_introduction_2007} or \cite{moller1995multivariate}. We leave this for future research.
\end{remark}
In what follows, we use the same function $h$ to define $\psi_\alpha$ and a set of homogeneous Lagrange polynomials. %To specify which scaling is used we will sometimes call $\l_1, \ldots, \l_\D$ a set of homogeneous Lagrange polynomials of degree $\alpha$ w.r.t.\ $h \in S_\alpha$. 
\begin{proposition} \label{prop:inj}
Let $\alpha \in \Cl(X)$ be such that no $\z_j$ is a basepoint of $S_\alpha$. Then
\begin{enumerate}
\item $\psi_\alpha$ is injective if and only if $I_\alpha = J_\alpha$. In this case $\HF_I(\alpha) \leq \D$,
\item $\psi_\alpha$ is surjective if and only if there exists a set of homogeneous Lagrange polynomials of degree $\alpha$. In this case $\HF_I(\alpha) \geq \D$.
\end{enumerate}
\end{proposition}
\begin{proof}
Let $f, h \in S_\alpha$ such that $\z_j \notin \V_X(h), j = 1, \ldots, \D$. If $\psi_\alpha$ is injective, then $ f \in J_\alpha \Rightarrow \psi_\alpha(f + I_\alpha) = 0 \Rightarrow f \in I_\alpha$. So $J_\alpha \subset I_\alpha$ and the other inclusion is trivial. Conversely, if $I_\alpha = J_\alpha$, then $\psi_\alpha(f + I_\alpha) = 0 \Rightarrow f \in J_\alpha \Rightarrow f \in I_\alpha$, so $\psi_\alpha$ is injective. The corresponding statement about $\HF_I$ follows easily.\\
If $\psi_\alpha$ is surjective, take $\l_j \in \psi_{\alpha}^{-1}(e_j)$. Conversely, if $\l_j, j = 1, \ldots, \D$ is a set of homogeneous Lagrange polynomials of degree $\alpha$, $\psi_\alpha(\l_j + I_\alpha) = e_j$ and $\psi_\alpha$ is surjective. Again, the statement about $\HF_I$ follows easily.
\end{proof}

\begin{corollary} \label{cor:radical}
If $\alpha \in \Pic(X)$ is ample\footnote{A divisor $D$ and its degree $\alpha = [D]$ are called \textit{very ample} if  $D$ is basepoint free and $X \rightarrow \PP(\Gamma(X, \OO_X(D))^\vee)$ is a closed embedding. If $kD$ (or $k\alpha$) is very ample for some $k \geq 1$, then $D$ (or $\alpha$) is called \textit{ample}. See \cite[Chapter 6]{cox2011toric} for definitions and properties.} and $I$ is radical, then $\psi_\alpha$ is injective. 
\end{corollary}
\begin{proof}
In this case $I = \II(\overline{G\cdot z_1} \cup \cdots \cup \overline{G\cdot z_\D} \cup Z')$ by the Nullstellensatz. Take $f \in J_\alpha$. Since any polynomial in $S_\alpha$ for $\alpha$ ample vanishes on $Z$ ($S_\alpha \subset K$, see e.g.\ \cite{soprounov2005toric}), $f$ vanishes on $Z' \subset Z$. Therefore $f \in I_\alpha$ and $J_\alpha \subset I_\alpha \subset J_\alpha$. Now apply Proposition \ref{prop:inj}.
\end{proof}
%We now give a geometric interpretation of what it means for a set of Lagrange polynomials to exist in $S_\alpha$. 
The following proposition shows that the existence of homogeneous Lagrange polynomials of degree $\alpha \in \Cl(X)$ is equivalent to the fact that the points $\Phi_\alpha(z_j)$ span a linear space of dimension $\D-1$ in $\PP^{n_\alpha -1}$. Let $p_j \in \C^{n_\alpha}$ be a set of homogeneous coordinates (in the standard sense) of $\Phi_\alpha(z_j) \in \PP^{n_\alpha - 1}$ and define the matrix $L_\alpha = [ p_1 ~ \cdots ~ p_\D ] \in \C^{n_{\alpha} \times \D}$.
\begin{proposition} \label{prop:Lalpha}
Let $\alpha \in \Cl(X)$ be such that no $\z_j$ is a basepoint of $S_\alpha$. There exists a set of Lagrange polynomials of degree $\alpha$ if and only if $L_\alpha$ has rank $\D$.
\end{proposition}
\begin{proof}
The rank of $L_\alpha$ is $\D$ if and only if there exists a left inverse matrix $L_{\alpha}^\dagger \in \C^{\D \times n_\alpha}$ such that $L_\alpha^\dagger L_\alpha = \id_\D$ is the $\D \times \D$ identity matrix. We will show that this is equivalent to the existence of a set of homogeneous Lagrange polynomials of degree $\alpha$. Suppose that $L_\alpha^\dagger$ exists. The rows of $L_\alpha^\dagger$ should be interpreted as elements of $S_\alpha$ represented in the basis $\{x^{b_1}, \ldots, x^{b_{n_\alpha}} \}$. The columns of $L_\alpha$ are elements of $S_\alpha^\vee$ represented in the dual basis. Let the $j$-th row of $L_{\alpha}^\dagger$ correspond to $\tilde{\l}_j \in S_\alpha$. It is clear from $L_\alpha^\dagger L_\alpha = \id_\D$ that 
$$\pair{\tilde{\l}_j, p_i} = \tilde{\l}_j(z_i) =  \begin{cases}
1 & i = j, \\
0 & \textup{otherwise.}
\end{cases}$$  
By Lemma \ref{lem:nonzero}, there is $h\in S_\alpha$ such that $h(z_j) \neq 0, j = 1, \ldots, \D$. Then $\l_j = h(z_j)\tilde{\l}_j, j = 1, \ldots, \D$ are a set of homogeneous Lagrange polynomials. Conversely, if a set of homogeneous Lagrange polynomials exists, construct a matrix $\tilde{L}_\alpha^\dagger$ by plugging the coefficients of $\l_j$ into the $j$-th row. Then there is $h \in S_\alpha$ such that $\tilde{L}_\alpha^\dagger L_\alpha = \diag(h(z_1), \ldots, h(z_\D))$ is an invertible diagonal matrix. The left inverse is $L_\alpha^\dagger =   \diag(h(z_1), \ldots, h(z_\D))^{-1} \tilde{L}_\alpha^\dagger$.
\end{proof}
%Proposition \ref{prop:Lalpha} can be generalized to the case where $\Phi_\alpha$ has basepoints, assuming it is well-defined at $z_1, \ldots, z_\D$. 
Based on these results, we make the following definition. 
\begin{definition}[Regularity]
The regularity $\Reg(I) \subset \Cl(X)$ of $I$ is the subset of degrees $\alpha \in \Cl(X)$ for which no $\z_j$ is a basepoint of $S_\alpha$ and the following equivalent conditions are satisfied: 
\begin{enumerate}
\item $\psi_\alpha$ is an isomorphism,
\item $\HF_I(\alpha) = \D$ and $I_\alpha = J_\alpha$, 
\item $\HF_I(\alpha) = \D$ and there exists a set of homogeneous Lagrange polynomials of degree $\alpha$, 
\item $I_\alpha = J_\alpha$ and there exists a set of homogeneous Lagrange polynomials of degree $\alpha$.
\end{enumerate}
\end{definition}

\begin{theorem}
If $\alpha \in \Reg(I)$, $\alpha_0 \in \Cl(X)_+$ is such that no $\z_j$ is a basepoint of $S_{\alpha_0}$ and $\HF_I(\alpha + \alpha_0) = \D$, then $\alpha + \alpha_0 \in \Reg(I)$. 
\end{theorem}
\begin{proof}
Let $\l_j, j = 1, \ldots, \D$ be a set of homogeneous Lagrange polynomials of degree $\alpha$ w.r.t.\ $h \in S_\alpha$. It is easy to verify that for generic $h_0 \in S_{\alpha_0}$, $h_0 \l_j, j = 1, \ldots, \D$ is a set of homogeneous Lagrange polynomials of degree $\alpha + \alpha_0$ w.r.t.\ $hh_0$.
\end{proof}

If $\alpha \in \Pic(X)$ is basepoint free and $\HF_I(\alpha) = \D$, then to show that $\alpha \in \Reg(I)$, by Proposition \ref{prop:Lalpha} it suffices to show that $L_\alpha$ is of rank $\D$. If $\alpha$ is `large enough' (the associated polytope has enough lattice points), this seems reasonable to expect. Alternatively, by Proposition \ref{prop:inj} it suffices to show that $I_\alpha = J_\alpha$. Based on experimental evidence we propose the following conjecture.
\begin{conjecture} \label{conj}
Let $I = \ideal{f_1, \ldots, f_n} \subset S$ be a homogeneous ideal obtained as in Section \ref{sec:setup} such that $V_X(I)$ is a zero-dimensional subscheme of $U \subset X$. Let $\alpha_i = \deg(f_i) \in \Pic(X)$ be the basepoint free degrees of the generators. Then 
$ \alpha_0 + \alpha_1 + \ldots + \alpha_n \in \Reg(I)$
for all $\alpha_0 \in \Cl(X)_+$ such that no $\z_j$ is a basepoint of $S_{\alpha_0}$.
\end{conjecture}
In the rest of this section, we prove some weaker results to support Conjecture \ref{conj} and we continue our running example by investigating the regularity.

We consider the question for which $\alpha \in \Cl(X)$ we have $\HF_I(\alpha) = \D$. The following theorem generalizes Theorem 3.16 in \cite{csahin2016multigraded} in the case where $Z$ is small enough.
\begin{theorem} \label{thm:koszul1}
Let $I = \ideal{f_1, \ldots, f_n} \subset S$ be a homogeneous ideal obtained as in Section \ref{sec:setup} such that $V_X(I)$ is a zero-dimensional subscheme of $U \subset X$. Let $\alpha_i = \deg(f_i) \in \Pic(X)$ be the basepoint free degrees of the generators. If 
$\codim(Z) \geq n$
then for all basepoint free $\alpha_0 \in \Pic(X)_+$, $\HF_I(\alpha_0 + \alpha_1 + \ldots + \alpha_n) = \D$. 
\end{theorem}
\begin{proof}
Consider the Koszul complex 
$$ 0 \rightarrow S(- \sum_{i=1}^n \alpha_i) \rightarrow \bigoplus_{\substack{\I \subset \{1,\ldots,n\}\\ |\I| = n - 1}} S(-\alpha_\I) \rightarrow \cdots \rightarrow \bigoplus_{\substack{\I \subset \{1,\ldots,n\}\\ |\I| = 2}} S(-\alpha_\I) \rightarrow \bigoplus_{i=1}^n S(-\alpha_i) \rightarrow S $$
where $\alpha_\I = \sum_{i \in \I} \alpha_i$ and $S(-\alpha)$ is the Cox ring with twisted grading: $S(-\alpha)_\beta = S(\beta - \alpha)$. Since the orbit closures $\overline{G \cdot z_j}$ have dimension $k-n$ and by assumption $\dim(Z) \leq k-n$, the $f_i$ form a regular sequence in $S$. Hence the Koszul complex is exact. Restricting to the degree $\alpha = \alpha_0 + \alpha_1 + \ldots + \alpha_n$ part we get 
$$ 0 \rightarrow S(\alpha_0) \rightarrow \bigoplus_{\substack{\I \subset \{1,\ldots,n\}\\ |\I| = n - 1}} S( \alpha -\alpha_\I) \rightarrow \cdots \rightarrow \bigoplus_{\substack{\I \subset \{1,\ldots,n\}\\ |\I| = 2}} S(\alpha -\alpha_\I) \rightarrow \bigoplus_{i=1}^n S(\alpha -\alpha_i) \rightarrow S_\alpha. $$
Since $\alpha_0$ is basepoint free, it corresponds to a polytope $P_0$ and we have by \eqref{eq:gradedpiece} and \eqref{eq:sumofdiv}
$$ \dim_\C(S_{\alpha_0 + \alpha_\I}) = |(P_0+P_{\I}) \cap M|$$ 
with $P_\I = \sum_{i \in \I} P_i$ for any subset $\I \subset \{0, \ldots, n\}$. Counting dimensions we get 
$$ \dim_\C((S/I)_\alpha) = \sum_{\ell = 0}^n (-1)^{n-\ell} \sum_{\substack{\I \subset \{1,\ldots,n\}\\ |\I| = \ell}} |(P_0 + P_\I) \cap M|,$$
and the right hand side is the formula \eqref{eq:mvpoints} for the mixed volume $\D = \MV(P_1, \ldots, P_n)$. 
\end{proof}
Note that the conditions of Theorem \ref{thm:koszul1} are satisfied by all toric surfaces ($n = 2$). Here is an analogous result for the case where the system is `unmixed' (in some sense) and the corresponding polytope is normal. 
\begin{theorem} \label{thm:koszul2}
Let $I = \ideal{f_1, \ldots, f_n} \subset S$ be a homogeneous ideal obtained as in Section \ref{sec:setup} such that $V_X(I)$ is a zero-dimensional subscheme of $X$. Let $\alpha_i = \deg(f_i) \in \Pic(X)$ be the basepoint free degrees of the generators. If there is a basepoint free degree $\alpha_\star \in \Pic(X)$ corresponding to a normal polytope, such that $\alpha_i = t_i \alpha_\star$ for positive integers $t_i$, then $\HF_I(t \alpha_\star) = \D$ for $t \geq \sum_{i=1}^n t_i$.
\end{theorem}
\begin{proof}
The assumption on $\alpha_i$ implies that $P_i = t_i P_\star+ m_i$ for a normal polytope $P_\star$, lattice points $m_i$ and positive integers $t_i$. We can assume without loss of generality that $m_i = 0, i = 1, \ldots, n$. We consider the embedding $X_\A \subset \PP^{|\A|-1}$ of $X$ where $\A = P_\star \cap M$. More precisely, $X_\A$ is the image of $\Phi_{\alpha_\star}$ \cite[Proposition 5.4.7]{cox2011toric}. Let $u_m, m \in \A$ be homogeneous coordinates on $\PP^{n_{\alpha_\star}-1} = \PP^{|\A| -1}$. The toric ideal of $X_\A$ is denoted $I_\A \subset \C[u_m, m \in \A]$ and the $\Z$-graded coordinate ring of $X_\A$ is $\C[X_\A] = \C[u_m, m \in \A]/I_\A$. %The $f_i$ define hypersurfaces of degree $t_i$ in $\PP^{|\A| - 1}$ given by $h_i = 0$. 
By \cite[Theorem 5.4.8]{cox2011toric}, we have $S_{\alpha_i} \simeq \C[X_\A]_{t_i}$ and $f_i \in S_{\alpha_i}$ corresponds to an element $h_i + I_\A \in \C[X_\A]_{t_i}$.
We define the homogeneous ideal $ I' = \ideal{h_1 + I_\A,\ldots, h_n + I_\A} \subset \C[X_\A]$. By assumption, $I'$ defines a 0-dimensional subscheme of $X_\A$, so $h_1 + I_\A, \ldots, h_n + I_\A$ is a regular sequence in $\C[X_\A]$. The ring $\C[X_\A]$ is arithmetically Cohen-Macaulay \cite[Exercise 9.2.8]{cox2011toric}, so the corresponding Koszul complex 
$$ 0 \rightarrow K_n \rightarrow K_{n-1} \rightarrow \cdots \rightarrow K_2 \rightarrow K_1 \rightarrow \C[X_\A] \quad \textup{ with } \quad K_t = \bigoplus_{\substack{\I \subset \{1,\ldots,n\}\\ |\I| = t}} \C[X_\A](-\sum_{i \in \I} t_i)$$
is exact. Since $P_\star$ is a normal polytope, we have $\dim_\C(\C[X_\A]_t) = |tP_\star \cap M|$. Counting dimensions and using the same formula as before for $\D = \MV(P_1, \ldots, P_n) = \MV(P_\star, \ldots, P_\star)$ we find that $\dim_\C((\C[X_\A]/I')_t) = \D$ for $t \geq \sum_{i=1}^n t_i$. Combining this with $(\C[X_\A]/I')_t \simeq (S/I)_{t \alpha_\star}$ (see \cite[Theorem 5.4.8]{cox2011toric}) we get the desired result.
\end{proof}
We note that in the case where $X$ is a product of projective spaces, stronger bounds than those of Theorem \ref{thm:koszul1} and Theorem \ref{thm:koszul2} are known \cite{bender2018towards}. 
\begin{example}
We continue Example \ref{ex:hirz2}. The polytope $P = P_1 + P_2$ (shown in Figure \ref{fig:hirzebruch2}) has 12 lattice points. Therefore $n_\alpha = 12$, with $\alpha = [D_P] \in \Pic(X)$. Since $\D = 3$, $L_\alpha$ is a $12 \times 3$ matrix. Its rows are indexed by the monomials spanning $S_\alpha$, and its columns by the representatives $z_j$. The transpose is given by 

\begin{small}
$$
L_\alpha^\top = ~ ~
\begin{blockarray}{ccccccccccccc}
\rot{$x_3x_4^2$} &
\rot{$x_1x_4^2$} & 
\rot{$x_2x_3^3x_4$} & 
\rot{$x_1x_2x_3^2x_4$} & 
\rot{$x_1^2x_2x_3x_4$} &
\rot{$x_1^3x_2x_4$} & 
\rot{$x_2^2x_3^5$} & 
\rot{$x_1x_2^2x_3^4$} & 
\rot{$x_1^2x_2^2x_3^3$} & 
\rot{$x_1^3x_2^2x_3^2$} & 
\rot{$x_1^4x_2^2x_3$} & 
\rot{$x_1^5x_2^2$}\\
\begin{block}{[cccccccccccc]c}
  ~1~ & -1~ & -1~ & 1~ & -1~ & 1~ & 1~ & -1~ & 1~ & -1~ & 1~ & -1~ &z_1~\\
  ~1~ & 0~ & -1~ & 0~ & 0~ & 0~ & 1~ & 0~ & 0~ & 0~ & 0~ & 0~ &z_2~\\
  ~0~& 1~& 0~&0~&0~&-1~&0~&0~&0~&0~&0~&1~&z_3~\\
\end{block}
\end{blockarray}.
$$
\end{small}
\noindent Consider $h = 39(x_3x_4^2 - x_1x_4^2) \in S_\alpha$ and note that $h(z_j) \neq 0$ for all $j$. A set of homogeneous Lagrange polynomials w.r.t.\ $h$ is given by

\begin{small}
$$
\frac{2~\tilde{L}_\alpha^\dagger}{13} = ~ ~
\begin{blockarray}{ccccccccccccc}
\rot{$x_3x_4^2$} &
\rot{$x_1x_4^2$} & 
\rot{$x_2x_3^3x_4$} & 
\rot{$x_1x_2x_3^2x_4$} & 
\rot{$x_1^2x_2x_3x_4$} &
\rot{$x_1^3x_2x_4$} & 
\rot{$x_2^2x_3^5$} & 
\rot{$x_1x_2^2x_3^4$} & 
\rot{$x_1^2x_2^2x_3^3$} & 
\rot{$x_1^3x_2^2x_3^2$} & 
\rot{$x_1^4x_2^2x_3$} & 
\rot{$x_1^5x_2^2$}\\
\begin{block}{[cccccccccccc]c}
 ~0~&0~&0~&2~&-2~&0~&0~&-2~&2~&-2~&2~&0 ~&\ell_1\\
 ~2~&0~&-2~&-1~&1~&0~&2~&1~&-1~&1~&-1~&0 ~&\ell_2\\
 ~0~&2~&0~&1~&-1~&-2~&0~&-1~&1~&-1~&1~&2 ~&\ell_3\\
\end{block}
\end{blockarray},
$$
\end{small}
\noindent which is related to the pseudo inverse of $L_\alpha$ by $$L_\alpha^\dagger = \diag(h(z_1), h(z_2), h(z_3))^{-1} \tilde{L}^\dagger_\alpha = \diag(1/78,1/39,1/39) \tilde{L}^\dagger_\alpha.$$
 To check that $I_\alpha = J_\alpha$ we compute $\HF_I(\alpha) = \HF_J(\alpha) = 3$. Hence we have $\alpha \in \Reg(I)$. In fact, in this example $I$ is radical and $\alpha$ is ample, so $I_\alpha = J_\alpha$ follows from Corollary \ref{cor:radical}. 
\end{example}

\section{A toric eigenvalue, eigenvector theorem} \label{sec:toriceval}
In this section, we will work with multiplication maps between graded pieces of $S/I$. Again, $I$ is a homogeneous ideal in $S$ obtained as in Section \ref{sec:setup} satisfying Assumptions \ref{ass:1}-\ref{ass:3}. For $\alpha, \alpha_0 \in \Cl(X)_+$, a homogeneous element $g \in S_{\alpha_0}$ defines a linear map 
$$ M_g : (S/I)_\alpha \rightarrow (S/I)_{\alpha + \alpha_0} : f + I_\alpha \mapsto gf + I_{\alpha + \alpha_0}$$
representing `multiplication with $g$'. Just as in the affine case, these multiplication maps will be the key ingredient to formulate our root finding problem as a linear algebra problem. We state a toric version of the eigenvalue, eigenvector theorem and show how the eigenvalues can be used to recover homogeneous coordinates of the solutions and equations for the corresponding $G$-orbits. Our main result uses the following Lemma.
\begin{lemma} \label{lem:isom}
Let $\alpha, \alpha_0 \in \Cl(X)_+$ be such that $\alpha, \alpha + \alpha_0 \in \Reg(I)$ and no $\z_j$ is a basepoint of $S_{\alpha_0}$. Then for generic $h_0 \in S_{\alpha_0}$, $\m_{h_0}: (S/I)_\alpha \rightarrow (S/I)_{\alpha + \alpha_0} : f + I_\alpha \mapsto h_0f + I_{\alpha + \alpha_0}$ is an isomorphism of vector spaces. 
\end{lemma}
\begin{proof}
Let $\psi_\alpha$ be given as in \eqref{eq:evalmap} for some $h \in S_\alpha$. We can take $hh_0 \in S_{\alpha + \alpha_0}$ to define $\psi_{\alpha + \alpha_0}$. Then $\psi_{\alpha + \alpha_0} \circ \m_{h_0} = \psi_{\alpha}$ shows that $\m_{h_0}$ is invertible. 
\end{proof}
\begin{theorem}[Toric eigenvalue, eigenvector theorem] \label{thm:multiplication}
Let $\alpha, \alpha_0 \in \Cl(X)_+$ be such that $\alpha, \alpha+ \alpha_0 \in \Reg(I)$ and no $\z_j$ is a basepoint of $S_{\alpha_0}$. Then for any $g \in S_{\alpha_0}$, $\m_{g/h_0} = \m_{h_0}^{-1} \circ \m_g: (S/I)_\alpha \rightarrow (S/I)_{\alpha}$ has eigenpairs 
$$ \left( \frac{g}{h_0}(z_j), \l_j + I_\alpha \right), \qquad \left( v_j, \frac{g}{h_0}(z_j) \right), \quad j = 1,  \ldots, \D,$$
where the $\l_j + I_\alpha$ are cosets of homogeneous Lagrange polynomials of degree $\alpha$ and the $v_j$ are the dual basis of $(S/I)_\alpha^\vee$.
\end{theorem}
\begin{proof}
The map $\m_{h_0}$ is an isomorphism by Lemma \ref{lem:isom}. We define $\psi_\alpha$, $\psi_{\alpha + \alpha_0}$ as in \eqref{eq:evalmap} with $h \in S_\alpha$, $h h_0 \in S_{\alpha + \alpha_0}$ respectively. A straightforward computation shows that
$\psi_{\alpha+\alpha_0} \circ \m_{h_0} (\l_j + I_\alpha) = e_j$.
Analogously, we have $ \psi_{\alpha+\alpha_0} \circ \m_g (\l_j + I_\alpha) = \frac{g}{h_0}(z_j) e_j$.
It follows that $h_0(z_j) \m_g(\l_j + I_\alpha) = g(z_j) \m_{h_0}(\l_j + I_\alpha)$, and therefore 
$$ \m_{g/h_0}(\l_j + I_\alpha) = \frac{g}{h_0}(z_j) (\l_j+ I_\alpha),$$
which proves the statement about the right eigenpairs, since the $\l_j + I_\alpha$ are linearly independent. For the statement about the left eigenpairs, note that for any $f \in S_\alpha$
$$ v_j \circ \m_{g/h_0} (f + I_\alpha) = v_j \circ \m_{h_0}^{-1} (gf + I_{\alpha+\alpha_0})$$
and since $\m_{h_0}$ is an isomorphism, there is $\tilde{f} \in S_\alpha$ such that $gf - h_0 \tilde{f} \in I_{\alpha+\alpha_0}$. Therefore, for each $z_j \in \V(I)$ we have 
$$ \frac{gf - h_0 \tilde{f}}{h_0 h}(z_j) = 0 \Rightarrow \frac{\tilde{f}}{h}(z_j) = \frac{g}{h_0}(z_j) \frac{f}{h}(z_j)$$
and thus, since $\m_{h_0}^{-1}(gf+I_{\alpha+\alpha_0}) = \tilde{f} + I_\alpha$, we have 
$$v_j \circ \m_{g/h_0} (f + I_\alpha) = v_j (\tilde{f} + I_\alpha) = \frac{g}{h_0}(z_j) v_j(f + I_\alpha).$$
The $v_j$ are linearly independent, so this concludes the proof.
\end{proof}
Let $S_{\alpha_0} = \bigoplus_{i=1}^{n_{\alpha_0}} \C \cdot x^{b_i}$ where $\alpha_0 \in \Cl(X)_+$ is such that no $\z_j$ is a basepoint of $S_{\alpha_0}$. We now show how the eigenvalues of the $\m_{x^{b_i}/h_0}$ lead directly to a set of defining equations of $G \cdot z_j, j = 1, \ldots, \D$ if $\alpha_0$ is `large enough'. For every cone $\sigma \in \Sigma_P$, we define $U_\sigma = \C^k \setminus \V(x^{\hat{\sigma}}) = \Spec(S_{x^{\hat{\sigma}}})$. Note that $ \C^k \setminus Z = \bigcup_{\sigma \in \Sigma_P} U_\sigma$. Let $D_{\alpha_0}$ be a representative divisor: $\alpha_0 = [D_{\alpha_0}] =[ \sum_{i=1}^k a_{0,i} D_i]$. Let $P_0 \subset M_\R$ be the polytope $\{m \in M_\R ~|~ F^\top m + a_0 \geq 0 \}$. If $\alpha_0 \in \Pic(X)$, then for every $\sigma \in \Sigma_P$ there is $m_\sigma \in P_0 \cap M$ such that 
\begin{equation} \label{eq:distpoint}
\pair{u_i, m_\sigma} + a_{0,i} = 0, \quad \forall \rho_i \in \sigma(1),
\end{equation}
see for instance \cite[Lemma 3.4]{cox1995homogeneous} or \cite[Theorem 4.2.8]{cox2011toric}. If $D_{\alpha_0}$ is not Cartier, such an $m_\sigma$ does not exist for every cone $\sigma \in \Sigma_P$. We will denote the subset of cones for which $m_\sigma \in P_0$ satisfying \eqref{eq:distpoint} exists by $\widetilde{\Sigma}_P \subset \Sigma_P$. This set is nonempty since $\{0\} \in \widetilde{\Sigma}_P$. We write $P_0 \cap M = \{m_1, \ldots, m_{n_{\alpha_0}} \}$, $b_i = F^\top m_i + a_0$ and $b_\sigma = F^\top m_\sigma + a_0$. For all $\sigma \in \widetilde{\Sigma}_P$ we denote $P_0 \cap M - m_\sigma = \{m_1 - m_\sigma, \ldots, m_{n_{\alpha_0}} - m_\sigma \}$ (note that $0 \in P_0 \cap M - m_\sigma$) and 
$ \sigma^\vee = \{ m \in M_\R ~|~ \pair{u,m} \geq 0, \forall u \in \sigma \}, \sigma^\perp = \{ m \in M_\R ~|~ \pair{u,m} = 0, \forall u \in \sigma \}.$ We partition $P_0 \cap M - m_\sigma$ into 
$\M_\sigma^\perp = (P_0 \cap M - m_\sigma) \cap \sigma^\perp$
and $\M_\sigma = (P_0 \cap M - m_\sigma) \setminus \M_\sigma^\perp$. The inclusion
$$ \N \M_\sigma + \Z \M_\sigma^\perp = \left \{ \sum_{m \in \M_\sigma} c_m m + \sum_{m \in \M_\sigma^\perp} d_m m ~|~ c_m \in \N, d_m \in \Z \right \} \subset \sigma^\vee \cap M$$
is clear. In what follows, we will show that if equality holds for some simplicial $\sigma \in \widetilde{\Sigma}_P$ with $z_j \in U_\sigma$, then $\alpha_0$ is `large enough' to recover equations for $G \cdot z_j$ from the eigenvalues of the $M_{x^{b_i}/h_0}$.
\begin{theorem} \label{thm:orbiteq}
Let $z \in U_\sigma$ for a simplicial cone $\sigma \in \widetilde{\Sigma}_P$ such that $\pi(z)$ is not a basepoint of $S_{\alpha_0}$. Take $h_0 \in S_{\alpha_0}$ such that $h_0(z) \neq 0$ and let $\lambda_{i} = z^{b_i}/h_0(z), i = 1, \ldots, n_{\alpha_0}$. If $ \sigma^\vee \cap M = \N \M_\sigma + \Z \M_\sigma^\perp$, then $G \cdot z \subset U_\sigma$ is the subvariety defined by the ideal 
$$ \ideal{x^{b_i-b_\sigma} - \lambda_{i} \frac{h_0(x)}{x^{b_\sigma}}~|~ i = 1, \ldots, n_{\alpha_0}} \subset S_{x^{\hat{\sigma}}}.$$ 
\end{theorem}
To prove Theorem \ref{thm:orbiteq}, we need the following auxiliary lemma.
\begin{lemma} \label{lem:orbiteq}
Let $\sigma \in \widetilde{\Sigma}_P$ be a simplicial cone. For any point $z \in U_\sigma$, the orbit $G \cdot z$ is the subvariety defined by 
$$ G \cdot z = \{ x \in U_\sigma ~|~ x^{F^\top m} - z^{F^\top m}, m \in \sigma^\vee \cap M \} \subset U_\sigma.$$
If $\sigma^\vee \cap M = \N \{m_1, \ldots, m_\kappa \} + \Z \{m_{\kappa + 1}, \ldots, m_s \}$, then 
$$ \{ x \in U_\sigma ~|~ x^{F^\top m} - z^{F^\top m} , m \in \sigma^\vee \cap M \}= \{ x \in U_\sigma ~|~  x^{F^\top m_i} - z^{F^\top m_i}, i = 1, \ldots, s \} .$$
\end{lemma}
\begin{proof}
%Note that $m_{\kappa + 1}, \ldots, m_s$ are a $\Z$-basis for $\sigma^\perp \cap M$. 
Note that $x^{F^\top m} - z^{F^\top m} \in S_{x^{\hat{\sigma}}}, \forall m \in \sigma^\vee \cap M$ and $m_{\kappa+1}, \ldots, m_s \in \sigma^\perp \cap M$. The first statement is shown in the proof of Theorem 2.1 in \cite{cox1995homogeneous}. For the second statement, the inclusion `$\subset$' is obvious. To show the opposite inclusion, take $m \in \sigma^\vee \cap M$ and write $m = c_1m_1 + \ldots + c_s m_s$ with $c_1, \ldots, c_\kappa \in \N$, $c_{\kappa+1}, \ldots, c_s \in \Z$. Then 
$$ x^{F^\top m} = \prod_{i=1}^{\kappa} (x^{F^\top m_i})^{c_i} \prod_{j=\kappa + 1}^{s} (x^{F^\top m_j})^{c_j}$$
and if $x^{F^\top m_i} = z^{F^\top m_i}, i = 1, \ldots, s$, it follows that $x^{F^\top m} = z^{F^\top m}$.
\end{proof}

\begin{proof}[Proof of Theorem \ref{thm:orbiteq}]
It follows from Lemma \ref{lem:orbiteq} that $G\cdot z$ is the variety of 
$$ \ideal{ x^{F^\top (m_i - m_\sigma)} - z^{F^\top (m_i - m_\sigma)} ~|~ i = 1, \ldots, n_{\alpha_0}} = \ideal{ x^{b_i - b_\sigma} - z^{b_i - b_\sigma}~|~i = 1, \ldots, n_{\alpha_0}}.$$
Write $h_0(x) = \sum_{i=1}^{n_{\alpha_0}} c_i x^{b_i}, c_i \in \C$. It is easy to check that 
$$ \left ( \begin{bmatrix}
1 \\ & \ddots & \\ & & 1
\end{bmatrix}  - \begin{bmatrix}
\lambda_{1} \\ \vdots \\ \lambda_{n_{\alpha_0}}
\end{bmatrix} \begin{bmatrix}
c_1 & \hdots & c_{n_{\alpha_0}}
\end{bmatrix} \right) \begin{bmatrix}
x^{b_1 - b_\sigma} - z^{b_1 - b_\sigma} \\ \vdots \\  x^{b_{n_{\alpha_0}} - b_\sigma} - z^{b_{n_{\alpha_0}} - b_\sigma}
\end{bmatrix} = \begin{bmatrix}
x^{b_1 - b_\sigma} - \lambda_{1} \frac{h_0(x)}{x^{b_\sigma}} \\ \vdots \\  x^{b_{n_{\alpha_0}} -b_\sigma} - \lambda_{n_{\alpha_0}} \frac{h_0(x)}{x^{b_\sigma}}
\end{bmatrix}$$
and for generic $c_i$, the matrix on the left is invertible (it's invertible for $c_i = 0$, so the determinant is a nonzero polynomial in the $c_i$).  
\end{proof}

\begin{theorem} \label{thm:usefuleq}
Let $z \in U_\sigma$ with $\sigma \in \widetilde{\Sigma}_P$ simplicial be such that $\pi(z)$ is not a basepoint of $S_{\alpha_0}$ and $ \sigma^\vee \cap M = \N \M_\sigma + \Z \M_\sigma^\perp$. For generic $h_0 \in S_{\alpha_0}$ satisfying $h_0(z) \neq 0$, the variety 
$$Y_z = \V \left (x^{b_i} - \frac{z^{b_i}}{h_0(z)}, i = 1, \ldots, n_{\alpha_0} \right ) \subset \C^k$$
is nonempty and $Y_z \subset G \cdot z$.
\end{theorem}
The proof of Theorem \ref{thm:usefuleq} uses the following lemma. 
\begin{lemma} \label{lem:nontorsion}
If $\alpha_0 \in \Cl(X)_+$ is such that $ \sigma^\vee \cap M = \N \M_\sigma + \Z \M_\sigma^\perp$ for some $\sigma \in \widetilde{\Sigma}_P$, then $\alpha_0$ is not a torsion element of $\Cl(X)$.
\end{lemma}
\begin{proof}
Suppose $u \alpha_0 = 0$ for some $u >0$. Then $F^\top m + u a_0 = 0$ for some $m \in M$, and therefore $F^\top (m/u) + a_0 = 0$. Since $\Sigma_P$ is complete, this means that $P_0 = \{m/u\}$ and $P_0$ either has 1 lattice point (if $m/u \in M$, in which case $\alpha_0 = 0$), or it has none. Since $\alpha_0 \in \Cl(X)_+$, we can assume $0 \in P_0$ and this must be the only lattice point in $P_0$. Then $\sigma^\vee \cap M = \N \M_\sigma + \Z \M_\sigma^\perp = \{0\}$. But $\sigma^\vee$ has dimension $n$ because $\sigma$ is strongly convex (\cite[Proposition 1.2.12]{cox2011toric}), so this is a contradiction.
\end{proof}

\begin{proof}[Proof of Theorem \ref{thm:usefuleq}]
%We first show that $Y_z \subset U_\sigma$. %Suppose, without loss of generality, that the vertices $\textup{v}_1, \ldots, \textup{v}_\ell$ of the face $\F_\sigma \subset P_0$ corresponding to $\sigma$ correspond to the monomials $x^{b_1}, \ldots, x^{b_\ell} \in S_{\alpha_0}$. %Since $P_0$ is ample, every Cox variable $x_i$ corresponds to a facet of $P_0$. 
%Every face not containing $\F_\sigma$ does not contain one of the $\textup{v}_1, \ldots, \textup{v}_\ell$. Hence, every $x_i, \rho_i \notin \sigma(1)$ divides one of the $x^{b_i}, i = 1, \ldots, \ell$. It follows that 
%$ \left( \prod_{\rho_i \notin \sigma(1)} x_i \right )$ divides $\left (\prod_{i=1}^\ell x^{b_i} \right)$.
%For $x \in Y_z$, we have 
%$$ \prod_{i=1}^\ell x^{b_i} = \prod_{i=1}^\ell \frac{z^{b_i}}{h_0(z)},$$
%which is nonzero since $z \in U_\sigma$, so $x_i \neq 0$ for $\rho_i \notin \sigma(1)$ and $x \in U_\sigma$. 
%Let $\F_\sigma =  \{m \in P_0 ~|~ \pair{u_i,m}+a_{0,i} =0, \forall \rho_i \in \sigma(1) \}$. By Lemma \ref{lem:secondlem}, for every $\rho_i \notin \sigma(1)$, there is $m \in \F_\sigma \cap M$ such that $\pair{u_i,m} + a_{0,i} >0$. It follows that for every $x_i, \rho_i \notin \sigma(1)$ there is $x^{F^\top m + a_0} = x^{b}$ with $m \in \F_\sigma \cap M$ such that $x_i | x^{b}$. For $x \in Y_z$, we have $h_0(z) x^{b} = z^{b}$, with $z^{b} \neq 0$ since $z \in U_\sigma$. Therefore, for $x \in Y_z, \rho_i \notin \sigma(1)$ we have $x_i \neq 0$.
%Hence $Y_z \subset U_\sigma$. 
Since $\alpha_0$ is not a torsion element of $\Cl(X)$ (Lemma \ref{lem:nontorsion}), we have the exact sequence 
$$ 0 \longrightarrow \Z \longrightarrow \Cl(X) \longrightarrow \Cl(X)/(\Z \cdot \alpha_0) \longrightarrow 0$$
where $\Z \rightarrow \Cl(X)$ sends $u \mapsto u \alpha_0 \in \Cl(X)$. Taking $\Hom_\Z(-, \C^*)$ shows that $G \rightarrow \C^* : g \mapsto g^{\alpha_0}$ is surjective (because $\C^*$ is divisible). Therefore we can find $g \in G$ such that $g^{\alpha_0} = h_0(z)^{-1}$ and thus $h_0(g \cdot z) = 1$. Every $x \in Y_z$ satisfies $x^{b_i} - (g \cdot z)^{b_i} = 0, i = 1, \ldots, n_{\alpha_0}$: this follows from $(g \cdot z)^{b_i} = z^{b_i}/h_0(z)$. In particular, $x^{b_\sigma} = (g \cdot z)^{b_\sigma} \neq 0$ ($z \in U_\sigma$ and hence $g \cdot z \in U_\sigma$ since $U_\sigma$ is $G$-invariant) and therefore $x$ satisfies $x^{b_i - b_\sigma} = (g \cdot z)^{b_i - b_\sigma}, i = 1, \ldots n_{\alpha_0}$. By Lemma \ref{lem:orbiteq} it follows that $g \cdot z \in Y_z \subset G \cdot z$.
\end{proof}

Recall that we took $\alpha_0$ such that  no $\z_j$ is a basepoint of $S_{\alpha_0}$. We conclude this section with the following immediate corollary of Theorems \ref{thm:orbiteq} and \ref{thm:usefuleq}.
\begin{corollary} \label{cor:orbiteq}
Let $\lambda_{ij} = z_j^{b_i}/h_0(z_j)$ be the $j$-th eigenvalue of the $i$-th multiplication map $\m_{x^{b_i}/h_0}$. For $j = 1, \ldots, \D$, assume that $z_j \in U_{\sigma_j}$ for a simplicial cone $\sigma_j \in \widetilde{\Sigma}_P$ satisfying $ \sigma_j^\vee \cap M = \N \M_{\sigma_j} + \Z \M_{\sigma_j}^\perp$. The ideal 
$$ \ideal{x^{b_i-b_{\sigma_j}} - \lambda_{ij} \frac{h_0(x)}{x^{b_{\sigma_j}}}~|~ i = 1, \ldots, n_{\alpha_0}} \subset S_{x^{\hat{\sigma}_j}}$$ 
defines the orbit $G \cdot z_j \subset U_{\sigma_j}$, and for any point $z_j' \in \V(x^{b_i} - \lambda_{ij}, i = 1, \ldots, n_{\alpha_0}) \subset U_{\sigma_j}$, we have $\pi(z_j') = \z_j$. 
\end{corollary}
Corollary \ref{cor:orbiteq} implies that we can find homogeneous coordinates of the solutions from the eigenvalues $\lambda_{ij}$ by solving a system of binomial equations if $P_0$ `has enough lattice points'. Concretely, for every point $z_j$ there has to be a cone $\sigma_j \in \widetilde{\Sigma}_P$ such that $z_j \in U_{\sigma_j}$ and $ \sigma_j^\vee \cap M = \N \M_{\sigma_j} + \Z \M_{\sigma_j}^\perp$. Note that if all solutions are in the torus, then $z_j \in U_\sigma$ for $\sigma = \{0\} \in \widetilde{\Sigma}_P$ and this condition translates to the fact that $\Z(P_0 \cap M - m) = M$ for some $m \in P_0 \cap M$. If $P_0$ is very ample, then $\alpha_0 \in \Pic(X)$, so $\widetilde{\Sigma}_P = \Sigma_P$ and $ \sigma^\vee \cap M = \N \M_\sigma + \Z \M_\sigma^\perp$ holds for all $\sigma \in \Sigma_P$ \cite[Proposition 1.3.16]{cox2011toric}. We will elaborate on how to solve this system of binomial equations in the next section.

\section{Algorithm} \label{sec:alg}
In this section we present an eigenvalue algorithm for computing homogeneous coordinates of the points in $\V_X(I)$, where $I$ is an ideal satisfying Assumptions \ref{ass:1}-\ref{ass:3}. As in Theorem \ref{thm:multiplication}, let $\alpha, \alpha_0 \in \Cl(X)_+$ be such that $\alpha, \alpha + \alpha_0 \in \Reg(I)$ and no $\z_j$ is a basepoint of $S_{\alpha_0}$. In practice, we will take $\alpha = \alpha_1 + \cdots + \alpha_n$ where $\alpha_i = \deg(f_i)$ (by Conjecture \ref{conj}) and $\alpha_0$ `large enough' to recover all solutions (Corollary \ref{cor:orbiteq}). We denote 
$$ S_{\alpha_0} = \bigoplus_{i=1}^{n_{\alpha_0}} \C \cdot x^{b_i}.$$
We have that
$ \HF_I(\alpha) = \HF_I(\alpha + \alpha_0) = \D.$ Given a generic $h_0 \in S_{\alpha_0}$ and a surjective linear map $N: S_{\alpha+\alpha_0} \rightarrow \C^{\D}$ with $\ker N = I_{ \alpha+\alpha_0 }$, we define $$N_{h_0}: S_{\alpha} \rightarrow \C^{\D}: f \mapsto N(h_0f)$$ and assume that $N_{h_0}$ is surjective as well. Such a map $N$ can be computed directly from the input equations. We will come back to this later. Let $N^* : B \rightarrow \C^\D$ be the restriction of $N_{h_0}$ to a subspace $B \subset S_\alpha$ of dimension $\D$ such that $N^*$ is invertible, and let $$N_i: B \rightarrow \C^\D: f \mapsto N(x^{b_i} f), \qquad i = 1, \ldots, n_{\alpha_0}.$$
\begin{theorem} \label{thm:simtrans}
The map $\nu: B \simeq (S/I)_\alpha: g \mapsto g + I_\alpha$ is an isomorphism of vector spaces and the linear maps $(N^*)^{-1} \circ N_i: B \rightarrow B$ satisfy $\nu \circ (N^*)^{-1} \circ N_i = \m_{x^{b_i}/h_0} \circ \nu$ where $\m_{x^{b_i}/h_0}$ are the maps from Theorem \ref{thm:multiplication}.
\end{theorem}
\begin{proof}
By Lemma \ref{lem:isom}, $h_0 f \in I_{\alpha + \alpha_0}$ if and only if $f \in I_\alpha$. Therefore $\ker N_{h_0} = I_\alpha$. The first statement follows from $S_\alpha = B \oplus \ker N_{h_0}$. Since $\ker N = I_{\alpha+ \alpha_0}$, $N$ is well defined mod $I_{\alpha + \alpha_0}$. We define 
$$\tilde{N} : (S/I)_{\alpha + \alpha_0} \rightarrow \C^\D : f + I_{\alpha + \alpha_0} \mapsto N(f).$$ 
Since $\tilde{N}$ is a surjective linear map between $\D$-dimensional vector spaces, it is invertible.
For $g \in B$, $N^*(g) = N(h_0 g) = (\tilde{N} \circ \m_{h_0}) (g + I_\alpha)$ so $\nu \circ (N^*)^{-1} = (\tilde{N} \circ \m_{h_0})^{-1}$. Analogously we find $N_i(g) =  (\tilde{N} \circ \m_{x^{b_i}}) (g + I_\alpha)$. 
The theorem follows from
$ (\nu \circ (N^*)^{-1} \circ N_i )(g) = ((\tilde{N} \circ \m_{h_0})^{-1} \circ (\tilde{N} \circ \m_{x^{b_i}})) (g+ I_\alpha) = (\m_{h_0}^{-1} \circ \m_{x^{b_i}} \circ \nu)(g).$
\end{proof}
Theorem \ref{thm:simtrans} tells us that, identifying $B$ with $(S/I)_\alpha$, the homogeneous multiplication operators are given by $(N^*)^{-1} \circ N_i$. After fixing a basis $\B$ for $B$ the multiplication operators are commuting $\D \times \D$ matrices and we can compute their simultaneous diagonalization to find the values $\lambda_{ij} = z_j^{b_i}/h_0(z_j)$. 

We now show how the map $N$ can be computed from the input equations. Our strategy is based on techniques for computing Truncated Normal Forms (TNFs), as introduced in \cite{telen2017solving}. We use the notation 
$V = S_{\alpha + \alpha_0}$, $V_i = S_{\alpha + \alpha_0 - \alpha_i}$ and by the \textit{Resultant map} $\Res: V_1 \times \cdots \times V_n \rightarrow V$ we mean the linear map 
$$ (q_1, \ldots, q_n) \mapsto q_1f_1 + \ldots + q_nf_n. $$
When represented in matrix form, using monomial bases for the vector spaces involved, this map looks a lot like the resultant matrices coming from Macaulay and toric resultants \cite[Chapters 3 and 7]{cox2}.
Since $\im \Res = I_{\alpha+\alpha_0}$, the cokernel map of $\Res$ is a map $N :  V \rightarrow \C^\D \simeq V/\im \Res$ with the properties we need. 

The next step is to find the homogeneous coordinates of $\V_X(I)$ from the $\lambda_{ij}$. Suppose that $z_j \in U_{\sigma_j}$ for $\sigma_j \in \widetilde{\Sigma}_P$ and that $\alpha_0$ is such that $\sigma_j^\vee \cap M = \N \M_{\sigma_j} + \Z \M_{\sigma_j}^\perp$. By Corollary \ref{cor:orbiteq}, it remains to compute one point on the variety $\V(x^{b_i} - \lambda_{ij}, i = 1, \ldots, n_{\alpha_0} )$ for $j = 1, \ldots, \D$. If $\z_j \in T_X$, we can do this efficiently using only linear algebra as follows. Let $A = [b_1 ~ \cdots ~ b_{n_{\alpha_0}} ] \in \Z^{k \times n_{\alpha_0}}$ be the matrix of exponents and compute its Smith normal form: $U A V = S$ with $U,V$ unimodular and $S = [\diag(m_1, \ldots, m_r, 0, \ldots, 0) ~ 0 ] \in \Z^{k \times n_{\alpha_0}} $, where $m_i | m_{i+1}$. We make the substitution of variables $x_\ell = y_1^{U_{1 \ell}} \cdots y_k^{U_{k \ell}}$ to obtain the equivalent system of equations given by $y^{Ub_i} = \lambda_{ij}$. Applying the invertible transformation given by the matrix $V$, this simplifies to 
$$ y_\ell^{m_\ell} = \prod_{i=1}^{n_{\alpha_0}} \lambda_{ij}^{V_{i \ell}}, ~ \ell = 1, \ldots, r \quad \textup{  and  }  \quad 1 = \prod_{i=1}^{n_{\alpha_0}} \lambda_{ij}^{V_{i\ell}}, ~ r < \ell \leq k.$$
This imposes no conditions on $y_\ell, \ell > r$, so we can put $y_\ell = 1, \ell > r$. Taking the logarithm then shows that 
$$ \log y = [ \log y_1 ~ \cdots ~ \log y_k ] = \left [ w ~ 0_{k-r} \right ]$$ 
where $w = [\log \lambda_{1j} ~ \cdots \log \lambda_{n_{\alpha_0}j}]  [V_{:,1} ~ \cdots ~ V_{:,r} ]  \diag(1/m_1, \ldots, 1/m_r)$ and $0_{k-r}$ is a row vector of length $k-r$ with zero entries. To find the homogeneous coordinates, we only need to invert our change of coordinates and the logarithm: $$\log x = [\log x_1 ~ \cdots ~ \log x_k] = \log y ~ U, \qquad x_\ell = e^{\log x_\ell}, \ell = 1, \ldots, k .$$
Taking the logarithm has some advantages for the implementation: it reduces all computations to some matrix multiplications and it may prevent overflow. When $\z_j$ is not in the torus, some of the $\lambda_{ij}$ may be zero. In this case, to compute a point on 
$\V(x^{b_i} - \lambda_{ij}, i = 1, \ldots, n_{\alpha_0} )$, we may use a simple Newton iteration, for instance. In the \textit{nearly} degenerate situation, where $\lambda_{ij}$ is close to zero for some $i$, the approach above suffers from rounding errors. We take this into account by using the Smith normal form technique when $(\min_i |\lambda_{ij}|)/(\sum_{i=1}^{n_{\alpha_0}} |\lambda_{ij}|^2)^{1/2} > \tol$, where $|\cdot|$ denotes the modulus and $\tol$ is a predefined tolerance. This leads to Algorithm \ref{alg:coxcoords}.
\begin{algorithm}
\footnotesize
\caption{Computes the Cox coordinates of the points defined by $I = \ideal{f_1, \ldots, f_n}$}\label{alg:coxcoords}
\begin{algorithmic}[1]
\STATE{ $\Res \gets \textup{Matrix of the resultant map $V_1 \times \cdots \times V_n \rightarrow V$}$ }\label{constres}
\STATE{$N \gets \textup{Matrix of the cokernel $V \rightarrow \C^\D$ of $\Res$}$}
\STATE{$h_0 \gets \textup{Generic element of $S_{\alpha_0}$}$}
\STATE{Construct a matrix of $N_{h_0}$}
\STATE{Find $B \subset S_\alpha$ such that $(N_{h_0})_{|B}$ is invertible }\label{choiceofB}
\STATE{$N^* \gets (N_{h_0})_{|B}$}
\STATE{Construct a matrix of $N_i$, $1\leq i \leq n_{\alpha_0}$}
\FOR{$i=1,\ldots,n_{\alpha_0}$}
\STATE{ $\m_{x^{b_i}/h_0} \gets (N^*)^{-1} N_i$}
\ENDFOR
\STATE{ Compute $\lambda_{ij}, 1\leq i \leq n_{\alpha_0}, 1 \leq j \leq \D$ by sim.~diag.~of the $\m_{x^{b_i}/h_0}$}
\FOR{$j=1,\ldots,\D$}
\STATE{$\tilde{J}_j \gets \ideal{x^{b_i} - \lambda_{ij}, 1\leq i \leq n_{\alpha_0}} \subset S$}
\IF{$(\min_i |\lambda_{ij}|)/(\sum_{i=1}^{n_{\alpha_0}} |\lambda_{ij}|^2)^{1/2} > \tol$}
\STATE{Find one point $z_j \in \C^k$ on $\V(\tilde{J}_j)$ using SNF}
\ELSE{}
\STATE{Find one point $z_j \in \C^k$ on $\V(\tilde{J}_j)$ using Newton iteration}
\ENDIF
\ENDFOR
\STATE{\textbf{return} $z_1, \ldots, z_\D$}
\end{algorithmic}
\end{algorithm}

In line \ref{choiceofB} of the algorithm, the choice of the subspace $B$ is important for the numerical stability. A good choice is using QR factorization with optimal column pivoting as in \cite{telen2017stabilized, telen2017solving} which results in a basis for $(S/I)_\alpha$ consisting of monomials in $S$. An alternative is using the singular value decomposition, in which case $B$ is the orthogonal complement of $I_\alpha$ in $S_\alpha$ \cite[Section 3]{mourrain2018truncated}. We use the SVD for the experiments in this article. 

Algorithm \ref{alg:coxcoords} requires some computations involving polytopes. If one is interested in solving many systems with the same structure, it is advantageous to do these computations in an `offline' phase. The `online' algorithm then takes a basis of $S_{\alpha_0}$, $S_\alpha$ and $S_{\alpha + \alpha_0}$, a facet representation of $P$ and $P_0$ and the mixed volume $\D = \MV(P_1, \ldots, P_n)$ as inputs. The `offline' version of the algorithm computes all this information from the input equations, and generates an $\alpha_0$ such that $\Z(P_0 \cap M - m) = M$. This is enough to find (at least) all solutions in the torus by Corollary \ref{cor:orbiteq}.%The problem of automatically generating a `small' basepoint free degree $\alpha_0$ (small meaning that $P_0$ has few lattice points) is in general not so easy. A naive way of doing this is choosing $P_0 = P_j$ where $P_j$ is the polytope with the smallest number of lattice points among $P_1, \ldots, P_n$. This is often still to large. When $\Sigma$ is simplicial, every Weil divisor is $\Q$-Cartier and one can proceed by finding a small, full-dimensional polytope $\tilde{P}_0$ corresponding to a Weil divisor and then find the smallest integer $\ell$ such that the vertices of $P_0 = \ell \tilde{P}_0$ are in the lattice. 

To retrieve the coordinates in $(\C^*)^n$ of toric solutions from their homogeneous coordinates, we use the map \eqref{eq:GCQ}.

\begin{remark} \label{rem:complexity}
We conclude this Section with a remark on the complexity of Algorithm \ref{alg:coxcoords} as compared to the TNF algorithm of \cite{telen2017solving}. The first step in both algorithms is to compute the cokernel of a resultant map $\Res$. Since for both algorithms the monomials indexing the vector space $V$ in the definition of $\Res$ are the lattice points contained in a slightly enlarged version of the polytope $P = P_1 + \ldots + P_n$, this step takes roughly the same computation time for both algorithms. Even though the Cox ring has dimension $k > n$, the dimensions of its graded pieces correspond to the lattice points contained in $n$-dimensional polytopes. This is an important observation, because for larger problems, the computation of the cokernel of $\Res$ is the most expensive step of the algorithm. Next, both algorithms compute the multiplication matrices from this cokernel. This is more expensive for the algorithm in this paper: there are more multiplication maps. Another important difference is that for the TNF algorithm, the eigenvalues of the multiplication maps immediately give the coordinates of the solutions, whereas Algorithm \ref{alg:coxcoords} processes these eigenvalues to find the homogeneous coordinates (line 12-19). We conclude that Algorithm \ref{alg:coxcoords} is computationally more expensive overall. This should be considered the price that is payed for being more robust in nearly degenerate situations, which is the main objective in this paper. However, the increase of complexity is not dramatic: systems with thousands of solutions can be solved within reasonable time (see Subsection \ref{subsec:generic}), and there is certainly room for performance optimization in the current Matlab implementation, which is tested in the next Section. 
\end{remark}

\section{Examples}
\label{sec:examples}
Algorithm \ref{alg:coxcoords} is implemented in Matlab. Polymake is used for computations involving polytopes \cite{polymake:FPSAC_2009}, except for the mixed volume, which is computed using PHCpack \cite{verschelde1999algorithm}. In this section, we test the implementation on several examples and compare the results with those of some other polynomial system solvers. All computations are done in double precision arithmetic on an 8 GB RAM machine with an intel Core 17-6820HQ CPU working at 2.70 GHz. To measure the quality of an approximate solution, we compute the \textit{residual} as defined in \cite[Section 7]{telen2017stabilized} as a measure for the relative backward error. In double precision arithmetic, a residual of order $10^{-16}$ is the best one can hope for. The goal of the experiments is to show that Algorithm \ref{alg:coxcoords} meets our objectives: it finds \textit{all} solutions with \textit{good accuracy} within reasonable time. In particular, it does so for (nearly) degenerate systems with solutions on or near the exceptional divisors of $X$ that cannot be solved by other state of the art solvers.
\subsection{Points on $\H_2$}
We finish our running example by using Algorithm \ref{alg:coxcoords} to compute homogeneous coordinates of the solutions of the system defined in Example \ref{ex:hirz1}. We use $\tol = 10^{-12}$, $\alpha = \alpha_1 + \alpha_2$. For $\alpha_0 = \alpha_2$, Algorithm \ref{alg:coxcoords} finds three solutions. %The solution corresponding to $z_3$ (lying on $D_3$) is computed with a residual of $\approx 10^{-9}$, whereas the first two solutions have a residual of order $10^{-16}$. Using $\alpha_0 = \alpha_1$,
All three residuals are of order $10^{-16}$. %We leave the investigation of the influence of the choice of $\alpha_0$ on the accuracy as future research.

\begin{figure}
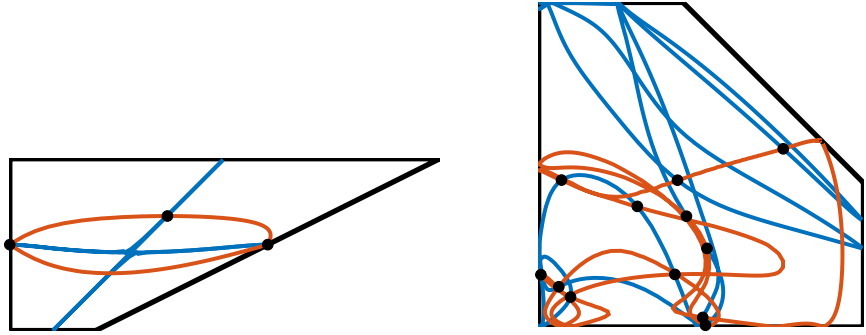

\centering
\input{moment_hirz.tex}
\hspace{1cm}
\input{chebmoment2.tex}
\caption{Left: images in $P$ of the real part of $\V(f_1)$ (\ref{bluecurve}) and $\V(f_2)$ (\ref{orangecurve}) from Example \ref{ex:hirz1} under the moment map $\mu$. The images of the computed real solutions are shown as black dots (\ref{sols}). Right: same picture for a different system.}
\label{fig:hirzmoment}
\end{figure}

To illustrate the results, we use the \textit{moment map} $$\mu: \C^k \setminus Z \rightarrow P : x \mapsto \frac{1}{\sum_{m \in P\cap M} |x^{F^\top m + a}|} \sum_{m \in P \cap M} |x^{F^\top m + a}| m,$$
where $|\cdot|$ denotes the modulus.
The map $\mu$ is constant on $G$-orbits and takes a point $x \in \C^k \setminus Z$ to a convex combination of the lattice points of $P$. It has the property that torus invariant prime divisors are sent to their corresponding facets and $(\C^*)^k$ is sent to the interior of $P$. More information can be found in \cite[Section 4.2]{fulton1993introduction} and \cite[Section 2]{sottile2017ibadan}. Figure \ref{fig:hirzmoment} shows that two of the computed solutions lie on divisors and one is in the torus. The image under $\mu$ of all of the solutions must lie on an intersection of the images of $\V(f_1) \setminus Z, \V(f_2) \setminus Z$ (but not all intersections correspond to solutions). As an illustration, we have included the same picture for a system with more solutions in the right part of the same figure. The polytopes for this system are $P_1 = [0,4] \times [0,4]$ and $P_2 = 5 \Delta_2$ where $\Delta_2$ is the standard simplex. There are $\D = 40$ solutions, 12 of them are real.

\subsection{A problem from computer vision} \label{subsec:vision}
\begin{figure}%
    \centering
    \subfloat[$P_1$]{{\includegraphics[width=2.2cm]{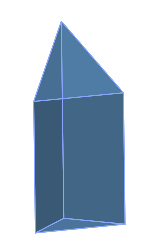} }}%
    \qquad
    \subfloat[$P_2$]{{\includegraphics[width=2.2cm]{P1.jpg} }}%
    \qquad
    \subfloat[$P_3$]{{\includegraphics[width=4cm]{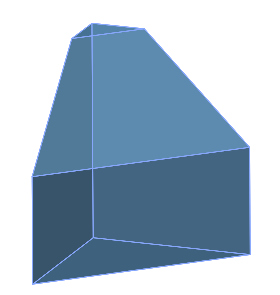} }}%
    \caption{Newton polytopes of the equations of the eight point radial distortion problem.}%
    \label{fig:polytopes8ptrad}%
\end{figure}

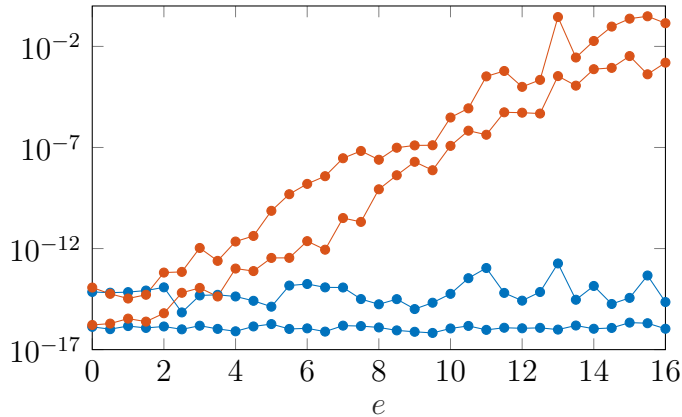
\begin{figure}[h!]
\centering
\input{residuals8pt.tex}
\caption{Minimal and maximal residual for different values of the parameter $e$ for the parametrized eight point radial distortion problem, for Algorithm \ref{alg:coxcoords} (\ref{blueres}) and the toric TNF algorithm (\ref{orangeres}).}
\label{fig:residuals8pt}
\end{figure}

One of the so-called `minimal problems' in computer vision is the problem of estimating radial distortion from eight point correspondences in two images. In \cite{kukelova2007minimal}, Kukelova and Pajdla propose a formulation of this problem as a system of 3 polynomial equations in 3 unknowns. The Newton polytopes are visualized in Figure \ref{fig:polytopes8ptrad}. The mixed volume is $\D = \MV(P_1,P_2,P_3) = 17$ and the matrix of facet normals is 
$$ F = \begin{bmatrix}
0 &-1& -1 &0& 1 &0 \\ 1 &-1 &-1 &0 &0 &0 \\ 0&0 &-1 &1 &0 &-1
\end{bmatrix},$$
so the Cox ring $S$ has dimension 6. We assign random real coefficients drawn from a standard normal distribution to all lattice points in the polytopes and solve the system using Algorithm \ref{alg:coxcoords}. We first run the offline version, which generates the polytope $P_0$. In this case, $P_0$ is the standard simplex. All 17 solutions are found with a residual of order $10^{-16}$ within $\pm 0.1$ s (using the online version of the algorithm). 
To show the robustness of Algorithm \ref{alg:coxcoords} in the nearly degenerate case, i.e.\ the case where there are solutions on or near the torus invariant prime divisors, we perform the following experiment. Consider the lattice points $$\F_3 = \{m \in P_1\cap M ~|~ \pair{u_3,m}+3 = 0 \}, \qquad \G_3 = (P_1 \cap M) \setminus \F_3.$$
The points in $\F_3$ are the lattice points on the facet of $P_1$ corresponding to $u_3 = (-1,-1,-1)$. Set 
$$ \hat{g}_i = \sum_{m \in \F_3} c_{m,i} \chi^m + \sum_{m \in \G_3} c_{m,i} \chi^m, \qquad i = 1,2$$
with $c_{m,i}$ real numbers drawn from a standard normal distribution. Now let $\f_1 = \hat{g}_1$ and $$\f_2(e) = \sum_{m \in \F_3} (10^{-e} c_{m,2} + (1-10^{-e})c_{m,1}) \chi^m + \sum_{m \in \G_3} c_{m,2} \chi^m, e \in [0, \infty).$$
The equation $\f_2 = 0$ is parametrized by the real parameter $e$. The third equation $\f_3 = 0$ is chosen randomly. When $e = 0$, $\f_2 = \hat{g}_2$ and the system is generic, as before. When $e \rightarrow \infty$, the part of $\f_2$ corresponding to $\F_3$ converges to the part of $\f_1$ corresponding to $\F_3$, meaning that there will be solutions `at infinity' on the divisor $D_3$. We solve the system for $e = 0, 1/2, 1, 3/2, \ldots, 16$ and compute both the maximal residual $r_{\max}$ and the minimal residual $r_{\min}$ for the 17 solutions found by Algorithm \ref{alg:coxcoords} with $\tol = 10^{-4}$ and the solutions found by the toric version of the Truncated Normal Form (TNF) algorithm \cite{telen2017solving}. The TNF solver computes the multiplication matrices for the input equations (in the classical sense) using heuristically `the best possible basis' from a numerical point of view. The numerical results in \cite{telen2017solving,mourrain2018truncated} motivate the choice of this method as a reference. The result of the experiment is shown in Figure \ref{fig:residuals8pt}. Note that not only the residuals of the solutions approaching the divisor deteriorate for the TNF algorithm.  Accuracy is lost on \textit{all} solutions. The reason is that even for the `best' basis selected by this algorithm, the computation of the classical multiplication matrices is ill-conditioned because the system is nearly degenerate. Looking at the computed Cox coordinates, we see that for three of the solutions, the coordinate $x_3$ goes to zero as $e$ increases, so 3 out of 17 solutions approach the divisor $D_3$. 

One can perform the same experiment for any other facet of $P_1$. However, in order to find the solutions on the divisors, the polytope $P_0$ must be large enough and it might not be sufficient that its lattice points generate the lattice (Corollary \ref{cor:orbiteq}). Repeating the same experiment, but this time using $\F_2$ instead of $\F_3$, the solutions in the torus are still found with good accuracy by Algorithm \ref{alg:coxcoords}. Accuracy is lost on the solutions approaching $D_2$. The reason is that the standard simplex does not `show' this facet. Using $P_0=\Conv((0,0,0),(1,0,0),(0,1,0),(0,0,1),(1,0,1),(0,1,1),(0,0,2))$ we find homogeneous coordinates of all solutions. 

\subsection{Generic problems} \label{subsec:generic}
To give an idea of the computation time and the type of systems Algorithm \ref{alg:coxcoords} can handle, we perform the following experiment. Consider the parameters $n$, NZ, $d_{\max} \in \N \setminus \{0\}$. For $j = 1,\ldots, n$ we generate a set $\A_j \subset \Z^n$ of NZ lattice points by selecting NZ points in $\N^n$ with coordinates drawn uniformly from $\{0,1,\ldots, d_{\max}\}$ and shifting these points by substracting the first point from all other points. Then for each $m \in \A_j$ we generate a random real number $c_{m,j}$ drawn from a standard normal distribution and we set 
$$ \f_j = \sum_{m \in \A_j} c_{m,j} \chi^m.$$
If two or more points $m \in \A_j$ coincide, we add the $c_{m,j}$ together, so NZ is an upper bound for the number of terms in $\f_j$. We use Algorithm \ref{alg:coxcoords} to compute the Cox coordinates of the solutions of the resulting system and their image under \eqref{eq:GCQ}. In Table \ref{tab:genericmixed} we report the number of solutions $\D$, the dimension $k$ of the Cox ring, the number $n_{\alpha_0}$ for the automatically generated $\alpha_0$, and, for both the offline and the online solver, the maximal residual $r_{\max}$, the geometric mean of the residuals of all solutions $r_{\mean}$ and the computation time $t$ (in seconds). The residuals are represented by $D_{\mean}= \ceil{-\log_{10} r_{\mean}}$ and $D_{\max} = \ceil{-\log_{10} r_{\max}}$.
\begin{table}[] 
\footnotesize
\centering
\begin{tabular}{ccc|ccc|ccc|ccc}
\multirow{2}{*}{$n$} & \multirow{2}{*}{NZ} & \multirow{2}{*}{$d_{\max}$} & \multirow{2}{*}{$\delta$} & \multirow{2}{*}{$k$} & \multirow{2}{*}{$n_{\alpha_0}$} & \multicolumn{3}{c}{OFFLINE}   & \multicolumn{3}{c}{ONLINE}    \\
                     &                     &                           &                           &                      &                        & $t$    & $D_{\mean}$ & $D_{\max}$ & $t$    & $D_{\mean}$ & $D_{\max}$ \\ \hline
2                    & 20                  & 10                        & 144                       & 12                   & 3                      & 1.9\text{e+}1 & 15        & 14       & 2.0\text{e-}1 & 15        & 14       \\
2                    & 20                  & 20                        & 505                       & 14                   & 4                      & 2.4\text{e+}1 & 14        & 12       & 1.9\text{e+}0 & 14        & 11       \\
2                    & 20                  & 30                        & 1268                      & 15                   & 3                      & 5.8\text{e+}1 & 14        & 12       & 1.9\text{e+}1 & 14        & 12       \\
2                    & 20                  & 40                        & 2390                      & 16                   & 3                      & 2.6\text{e+}2 & 14        & 11       & 1.4\text{e+}2 & 14        & 13       \\
2                    & 20                  & 50                        & 3275                      & 16                   & 3                      & 3.7\text{e+}2 & 14        & 12       & 2.3\text{e+}2 & 14        & 11       \\
2                    & 20                  & 60                        & 4469                      & 12                   & 3                      & 7.8\text{e+}2 & 11        & 7        & 5.2\text{e+}2 & 11        & 8        \\
2                    & 40                  & 30                        & 1522                      & 15                   & 3                      & 9.5\text{e+}1 & 14        & 11       & 3.4\text{e+}1 & 14        & 10       \\
2                    & 60                  & 30                        & 1670                      & 15                   & 4                      & 1.2\text{e+}2 & 14        & 12       & 5.3\text{e+}1 & 14        & 12       \\
2                    & 200                 & 30                        & 1672                      & 10                   & 3                      & 1.1\text{e+}2 & 15        & 10       & 6.0\text{e+}1 & 15        & 9        \\ \hline
3                    & 5                   & 3                         & 18                        & 21                   & 4                      & 2.2\text{e+}1 & 14        & 12       & 1.1\text{e-}1 & 15        & 13       \\
3                    & 5                   & 5                         & 136                       & 36                   & 4                      & 3.9\text{e+}1 & 14        & 9        & 6.3\text{e-}1 & 14        & 13       \\
3                    & 10                  & 5                         & 190                       & 60                   & 5                      & 3.5\text{e+}1 & 15        & 7        & 2.1\text{e+}0 & 15        & 11       \\
3                    & 10                  & 7                         & 592                       & 63                   & 5                      & 1.3\text{e+}2 & 14        & 10       & 3.2\text{e+}1 & 15        & 7        \\ \hline
4                    & 5                   & 3                         & 81                        & 106                  & 6                      & 6.9\text{e+}1 & 14        & 11       & 3.7\text{e+}1 & 14        & 11      
\end{tabular}
\caption{Results for generic systems with mixed supports.}
\label{tab:genericmixed}
\end{table}
It follows from Bernstein's second theorem \cite{bernstein,hustu} that solutions on divisors can only occur if the involved polytopes have common tropisms corresponding to positive dimensional faces. An important case in which this may happen is the unmixed case in which all input polytopes are equal. We repeat the experiment, but this time we keep the supports $\A = \A_1 = \ldots = \A_n$ fixed. Table \ref{tab:genericunmixed} shows some results. Of course, for this type of systems, the dimension of the Cox ring (or, equivalently, the number of facets of the Minkowski sum of the input polytopes) is lower and the system of binomial equations from Corollary \ref{cor:orbiteq} is easier to solve. %If $k$ is large, we observe that it gets trickier for the Newton iteration to converge to a point on $\V(J_j)$, so some solutions might be lost. A more sophisticated implementation of Newton's method could solve this problem. In these experiments, we did not use a factor $n-1$ in the determination of $D_{\rho_0}$ and observed numerically that this was not necessary. 
\begin{table}[h!]
\footnotesize
\centering
\begin{tabular}{ccc|ccc|ccc|ccc}
\multirow{2}{*}{$n$} & \multirow{2}{*}{NZ} & \multirow{2}{*}{$d_{\max}$} & \multirow{2}{*}{$\delta$} & \multirow{2}{*}{$k$} & \multirow{2}{*}{$n_{\alpha_0}$} & \multicolumn{3}{c}{OFFLINE}   & \multicolumn{3}{c}{ONLINE}    \\
                     &                     &                           &                           &                      &                        & $t$    & $D_{\mean}$ & $D_{\max}$ & $t$    & $D_{\mean}$ & $D_{\max}$ \\ \hline
2                    & 20                  & 60                        & 3638                      & 7                    & 3                      & 5.8\text{e+}2 & 13        & 11       & 3.8\text{e+}2 & 13        & 10       \\ \hline
3                    & 10                  & 10                        & 834                       & 14                   & 6                      & 3.5\text{e+}2 & 13        & 12       & 1.9\text{e+}2 & 13        & 12       \\ \hline
4                    & 6                   & 3                         & 15                        & 7                    & 8                      & 3.3\text{e+}1 & 15        & 15       & 8.4\text{e-}1 & 15        & 14       \\
4                    & 6                   & 4                         & 28                        & 6                    & 11                     & 4.3\text{e+}1 & 14        & 13       & 5.4\text{e+}0 & 15        & 14       \\
4                    & 6                   & 5                         & 216                       & 9                    & 7                      & 5.7\text{e+}2 & 12        & 11       & 2.7\text{e+}2 & 12        & 11       \\
4                    & 6                   & 6                         & 339                       & 8                    & 6                      & 1.5\text{e+}3 & 6         & 4        & 2.0\text{e+}3 & 6         & 5        \\ \hline
5                    & 6                   & 3                         & 10                        & 6                    & 8                      & 7.5\text{e+}1 & 15        & 14       & 1.0\text{e+}1 & 15        & 15      
\end{tabular}
\caption{Results for generic systems with unmixed supports.}
\label{tab:genericunmixed}
\end{table}

\subsection{Comparison with homotopy methods} \label{subsec:exphomotopy}
As discussed in the introduction, homotopy continuation methods provide very successful numerical solvers for systems of small degrees in large numbers of variables. Algebraic methods prove to be more robust in the case of high degrees and small dimensions, see for instance the numerical experiments in \cite{telen2017solving}. In this sense, these two important classes of numerical solvers are complementary to each other. As an illustration, we repeat the mixed experiment from Subsection \ref{subsec:generic} for three challenging 2-dimensional systems and compare the results with two homotopy implementations that are considered state of the art: Bertini (v1.6) \cite{bates2013numerically} and PHCpack (v2.4.64) \cite{verschelde1999algorithm}. For both these solvers, we use standard double precision settings and the backward errors of the computed solutions are of the order of the machine precision because these solvers intrinsically use Newton refinement. The results are reported in Table \ref{tab:homotopy}. For each solver, the number $\hat{\D}$ is the number of correctly computed solutions (with residual $<10^{-9}$). 
\begin{table}[]
\centering
\footnotesize
\begin{tabular}{ccc|c|ccccc|cc|cc}
$n$ & NZ  & $d_{\max}$ & $\D$    & \multicolumn{5}{c}{Algorithm \ref{alg:coxcoords}} & \multicolumn{2}{c}{PHCpack} & \multicolumn{2}{c}{Bertini} \\
  &     &      &      & $k$   & $n_{\alpha_0}$ & $t_{\textup{OFFLINE}}$   & $t_{\textup{ONLINE}}$    & $\hat{\D}$ & $t$           & $\hat{\D}$     & $t$             & $\hat{\D}$       \\ \hline
2 & 20  & 20   & 622  & 14  & 3  & 4.4e+1 & 2.8e+0 & 622  & 1.7e+0      & 597       & 2.2e+1       
& 605 \\
2 & 200 & 30   & 1700 & 14  & 3  & 1.5e+2 & 7.1e+1 & 1700 & 1.3e+1      & 1671      & 4.9e+2        & 1119        \\
2 & 800 & 40   & 3117 & 9   & 3  & 3.5e+2 & 2.3e+2 & 3117 & 7.7e+1      & 3055      &    7.6e+3            &   2832         
\end{tabular}
\caption{Results for generic systems using Algorithm \ref{alg:coxcoords} and the homotopy packages PHCpack and Bertini.}
\label{tab:homotopy}
\end{table}
Note that both homotopy solvers miss some solutions for all these problems. PHCpack is very efficient for this type of generic problems because it implements polyhedral homotopies \cite{hustu,verschelde1994homotopies}. This means in practice that exactly $\D$ paths are tracked. Bertini tracks 1258, 3135 and 6320 paths for the first, second and third problem respectively. This experiment shows that even for generic systems, for large $\D$ and small $n$ the state of the art homotopy algorithms do not find all solutions. The method introduced in this paper aims at solving (nearly) degenerate, non-generic systems. In practice, this often means that there are `large solutions'. To show that such situations cause trouble for homotopy methods, even for small $\D$, we consider the experiment of Subsection \ref{subsec:vision}. Solving the system for $e = 4.5$ using Algorithm \ref{alg:coxcoords} we find three solutions whose coordinates have a modulus of order $10^4$. PHCpack and Bertini both find only 14 solutions (the homotopy solvers give up on the paths converging to the `large solutions').
\section{Conclusion}
\label{sec:conclusion}
We have presented a toric eigenvalue, eigenvector theorem that allows to compute homogeneous coordinates of solutions of systems of Laurent polynomial equations (satisfying the assumptions in Section \ref{sec:setup}) on a natural toric compactification $X$ of $(\C^*)^n$. This results in a numerical linear algebra based algorithm that proves to be robust in the case of (nearly) degenerate systems with solutions on the torus invariant prime divisors. The algorithm is particularly successful for small dimensions $n$ and large degrees. It relies on a conjecture related to the regularity of $I$ (Conjecture \ref{conj}), which is checked numerically to be true in all of the presented experiments and supported by some weaker results in Section \ref{sec:lagreg}.

\section*{Acknowledgments}
I would like to thank David Cox for his many useful comments on an earlier version of this paper. I also want to thank Mateusz Michalek and Milena Wrobel for our discussions that led to Theorems \ref{thm:koszul1} and \ref{thm:koszul2} and their proofs, Tomas Pajdla and Zuzana Kukelova for suggesting the problem of Subsection \ref{subsec:vision}, Bernd Sturmfels for suggesting the title and Wouter Castryck, Alexander Lemmens, Bernard Mourrain, Marc Van Barel and Wim Veys for fruitful discussions. I am grateful to an anonymous referee for their detailed and useful suggestions.

%\section*{References}
\newpage
\footnotesize
\bibliography{references}
\bibliographystyle{alpha}

\end{document}

%% file: fan.tex
\definecolor{mycolor1}{rgb}{0.00000,0.44700,0.74100}%
\definecolor{mycolor2}{rgb}{0.85000,0.32500,0.09800}%
\definecolor{mycolor3}{rgb}{0.92900,0.69400,0.12500}%
\definecolor{mycolor4}{rgb}{0.4660,0.6740, 0.1880}%
		
		\begin{tikzpicture}[scale=0.8]
			\coordinate (O) at (0,0);
			\coordinate (u1) at (1,0);
			\coordinate (uu1) at (2.2,0);
			\coordinate (u12) at (2.2,-2.2);
			\coordinate (u2) at (0,-1);
			\coordinate (uu2) at (0,-2.2);
			\coordinate (u23) at (-2.2,-2.2);
			\coordinate (uu23) at (-2.2,2.2);
			\coordinate (u3) at (-1,2);
			\coordinate (uu3) at (-1.1,2.2);
			\coordinate (u4) at (0,1);
			\coordinate (uu4) at (0,2.2);
			\coordinate (u14) at (2.2,2.2);			
					\draw[-latex, very thick, red] (O) -> (u1); 
					\draw[-latex, very thick, red] (O) -> (u2);
					\draw[-latex, very thick, red] (O) -> (u3);
					\draw[-latex, very thick, red] (O) -> (u4); 
					\draw[opacity=0.2,fill = mycolor1,] (O)--(uu1)--(u12)--(uu2)--cycle;
					\draw[opacity=0.2,fill = mycolor2,] (O)--(uu2)--(u23)--(uu23)--(uu3)--cycle;
					\draw[opacity=0.2,fill = mycolor3,] (O)--(uu4)--(u14)--(uu1)--cycle;
					\draw[opacity=0.2,fill = mycolor4,] (O)--(uu3)--(uu4)--cycle;
					\node at (1.3,-0.3) {$u_1$};
					\node at (0.3,1.0) {$u_2$};
					\node at (-1.3,1.7) {$u_3$};
					\node at (-0.4,-0.8) {$u_4$};
					\node at (8,0) {$F = [u_1~u_2~u_3~u_4]= \begin{bmatrix}
					1&0&-1&0\\0&1&2&-1
					\end{bmatrix}$};
					\end{tikzpicture}

%% file: newtonpolygons.tex
\definecolor{mycolor1}{rgb}{0.00000,0.44700,0.74100}%
\definecolor{mycolor2}{rgb}{0.85000,0.32500,0.09800}%
\definecolor{mycolor3}{rgb}{0.92900,0.69400,0.12500}%

\begin{tikzpicture}
\begin{axis}[%
width=5.500in,
height=1.100in,
scale only axis,
xmin=-0.5,
xmax=12.5,
xtick = \empty,
ymin=-1,
ymax=2,
ytick = \empty,
axis background/.style={fill=white},
axis line style={draw=none} 
]
\addplot [color=mycolor1,solid,thick, fill opacity = 0.2, fill = mycolor1,forget plot]
  table[row sep=crcr]{%
0	0\\
1	0\\
3	1\\
0	1\\
0	0\\
};
\addplot[only marks,mark=*,mark size=1.5pt,mycolor1
        ]  coordinates {
    (0,0) (1,0) (3,1) (2,1) (1,1) (0,1)
};

\node at (axis cs:1.0,0.5) {$P_1$};
\node at (axis cs:3.5,0.5) {$+$};

\addplot [color=mycolor2,solid,thick, fill opacity = 0.2, fill = mycolor2,forget plot]
  table[row sep=crcr]{%
4	0\\
6	1\\
4	1\\
4	0\\
};
\addplot[only marks,mark=*,mark size=1.5pt,mycolor2
        ]  coordinates {
    (4,0) (6,1) (5,1) (4,1)
};
\node at (axis cs:4.5,0.6) {$P_2$};
\node at (axis cs:6.5,0.5) {$=$};

\addplot [color=mycolor3,solid,thick, fill opacity = 0.2, fill = mycolor3,forget plot]
  table[row sep=crcr]{%
7	-0.5\\
8	-0.5\\
12	1.5\\
7	1.5\\
7	-0.5\\
};

\addplot[only marks,mark=*,mark size=1.5pt,mycolor3
        ]  coordinates {
    (7,-0.5) (8,-0.5) (7,0.5) (8, 0.5) (9, 0.5) (10, 0.5) (7, 1.5) (8,1.5) (9,1.5) (10,1.5) (11,1.5) (12,1.5)
};

\node at (axis cs:8.4,0.6) {$P$};
\end{axis}
\end{tikzpicture}%

%% file: residuals8pt.tex
% This file was created by matlab2tikz.
%
%The latest updates can be retrieved from
%  http://www.mathworks.com/matlabcentral/fileexchange/22022-matlab2tikz-matlab2tikz
%where you can also make suggestions and rate matlab2tikz.
%
\definecolor{mycolor1}{rgb}{0.00000,0.44700,0.74100}%
\definecolor{mycolor2}{rgb}{0.85000,0.32500,0.09800}%
\begin{tikzpicture}

\begin{axis}[%
width=3in,
height=1.8in,
at={(0.772in,0.516in)},
scale only axis,
xmin=0,
xmax=16,
xlabel style={font=\color{white!15!black}},
xlabel={$e$},
ymode=log,
ymin=1e-17,
ymax=1,
yminorticks=true,
axis background/.style={fill=white}
]
\addplot [color=mycolor1, mark size=1.7pt, mark=*, mark options={solid, mycolor1}, forget plot]
  table[row sep=crcr]{%
0	7.11970312368113e-15\\
0.5	6.55116547609014e-15\\
1	6.99507996979959e-15\\
1.5	8.37245736709256e-15\\
2	1.20389513828192e-14\\
2.5	6.9183637171166e-16\\
3	4.79178997518004e-15\\
3.5	5.13838231404755e-15\\
4	4.26001368064262e-15\\
4.5	2.59681031877054e-15\\
5	1.33051498071503e-15\\
5.5	1.47710827415357e-14\\
6	1.76536478022386e-14\\
6.5	1.20092338265676e-14\\
7	1.18317818767741e-14\\
7.5	3.17449939150209e-15\\
8	1.74645169901965e-15\\
8.5	3.14374814003295e-15\\
9	1.01614265459418e-15\\
9.5	2.08943281511284e-15\\
10	5.74791855043334e-15\\
10.5	3.418871546494e-14\\
11	1.09039663299249e-13\\
11.5	6.38964560671947e-15\\
12	2.64282253273866e-15\\
12.5	7.08420082378827e-15\\
13	1.83497533241813e-13\\
13.5	2.93257421896477e-15\\
14	1.41225232492975e-14\\
14.5	1.81023028812561e-15\\
15	3.6295738431975e-15\\
15.5	4.64736692893318e-14\\
16	2.25862721261815e-15\\
}; \label{blueres}
\addplot [color=mycolor1, mark size=1.7pt, mark=*, mark options={solid, mycolor1}, forget plot]
  table[row sep=crcr]{%
0	1.33512936542337e-16\\
0.5	1.02932482579946e-16\\
1	1.46606429103363e-16\\
1.5	1.18473994624591e-16\\
2	1.39533098090834e-16\\
2.5	1.00707649391202e-16\\
3	1.5232996130677e-16\\
3.5	1.06160789552421e-16\\
4	8.07283330159285e-17\\
4.5	1.40240642221223e-16\\
5	1.83651737218573e-16\\
5.5	1.05395908165309e-16\\
6	1.12741258908529e-16\\
6.5	7.78228274680445e-17\\
7	1.53661608461273e-16\\
7.5	1.46369620625721e-16\\
8	1.25055751232726e-16\\
8.5	8.94155531536847e-17\\
9	7.6229192660466e-17\\
9.5	6.69553675311279e-17\\
10	1.11311277656623e-16\\
10.5	1.48914835184262e-16\\
11	9.48893776700998e-17\\
11.5	1.20794113972686e-16\\
12	1.13559528765453e-16\\
12.5	1.1941889472949e-16\\
13	9.7877992922458e-17\\
13.5	1.55817019067604e-16\\
14	1.07633042483581e-16\\
14.5	1.18107151554155e-16\\
15	2.16692008342442e-16\\
15.5	2.02474975292835e-16\\
16	1.07795580431296e-16\\
};
\addplot [color=mycolor2, mark size=1.7pt, mark=*, mark options={solid, mycolor2}, forget plot]
  table[row sep=crcr]{%
0	1.17692024462009e-14\\
0.5	5.74619061627348e-15\\
1	3.37254234522686e-15\\
1.5	5.17488978151407e-15\\
2	6.48009701276984e-14\\
2.5	7.00601985329477e-14\\
3	1.06973765249745e-12\\
3.5	2.46028420928425e-13\\
4	2.221115432497e-12\\
4.5	4.20641901915154e-12\\
5	7.34988713755362e-11\\
5.5	4.87897302636307e-10\\
6	1.58000358174065e-09\\
6.5	3.80150411570952e-09\\
7	2.9351650217648e-08\\
7.5	6.73309609843314e-08\\
8	2.47548925584109e-08\\
8.5	9.73533426320864e-08\\
9	1.26615481543515e-07\\
9.5	1.28797103477782e-07\\
10	3.03288135234367e-06\\
10.5	8.64431697379446e-06\\
11	0.000329038481978505\\
11.5	0.000609048881365017\\
12	9.9799007758833e-05\\
12.5	0.000220840067138713\\
13	0.283761868441605\\
13.5	0.00281998778475204\\
14	0.0185537520630187\\
14.5	0.0955639186123513\\
15	0.235896539108159\\
15.5	0.309203604028933\\
16	0.140556897292671\\
};
\addplot [color=mycolor2, mark size=1.7pt, mark=*, mark options={solid, mycolor2}, forget plot]
  table[row sep=crcr]{%
0	1.66001396696107e-16\\
0.5	1.9385210349436e-16\\
1	3.42425425231893e-16\\
1.5	2.42934338795001e-16\\
2	6.29431161562933e-16\\
2.5	6.47642941876884e-15\\
3	1.12683485935615e-14\\
3.5	4.24623661727479e-15\\
4	1.02334190604947e-13\\
4.5	7.68901925830485e-14\\
5	3.43758020219062e-13\\
5.5	3.46459840368101e-13\\
6	2.32761423355375e-12\\
6.5	8.77222545767569e-13\\
7	3.23330761066613e-11\\
7.5	2.0808432572322e-11\\
8	8.54501977654252e-10\\
8.5	4.22347199841747e-09\\
9	1.95563625228131e-08\\
9.5	7.48894466018289e-09\\
10	1.20346458268446e-07\\
10.5	6.77600293911804e-07\\
11	4.2716971553677e-07\\
11.5	5.49185695183361e-06\\
12	5.26325272397092e-06\\
12.5	4.76942289380576e-06\\
13	0.000340560356967094\\
13.5	0.000114589817591398\\
14	0.000747374368824812\\
14.5	0.000866357967373392\\
15	0.0033298021554058\\
15.5	0.000417262856207283\\
16	0.00156604746081053\\
}; \label{orangeres}
\end{axis}
\end{tikzpicture}%